\documentclass[preprint]{elsarticle}

\usepackage{mathrsfs}
\usepackage{amsfonts,amsmath,amssymb,amsthm}
\usepackage{graphicx}
\usepackage{latexsym}
\usepackage{tikz}
\usepackage{geometry}
\usepackage[mathscr]{eucal}
\usepackage[colorlinks=true]{hyperref}
\hypersetup{allcolors=[rgb]{0.1,0.1,0.4}}
\usepackage[all,cmtip]{xy}

\theoremstyle{plain}%
\newtheorem{theorem}{Theorem}
\newtheorem{proposition}[theorem]{Proposition}%
\newtheorem{lemma}[theorem]{Lemma}
\newtheorem{corollary}[theorem]{Corollary}

\newtheorem{remark}[theorem]{Remark}
\newtheorem{definition}[theorem]{Definition}

\newtheorem{example}[theorem]{Example}

\newcommand{\ua}{\mathord{\uparrow}}
\newcommand{\da}{\mathord{\downarrow}}

\newcommand{\ra}{\rightarrow }

\makeatletter
\def\@author#1{\g@addto@macro\elsauthors{\normalsize%
    \def\baselinestretch{1}%
    \upshape\authorsep#1\unskip\textsuperscript{%
      \ifx\@fnmark\@empty\else\unskip\sep\@fnmark\let\sep=,\fi
      \ifx\@corref\@empty\else\unskip\sep\@corref\let\sep=,\fi
      }%
    \def\authorsep{\unskip,\space}%
    \global\let\@fnmark\@empty
    \global\let\@corref\@empty  
    \global\let\sep\@empty}%
    \@eadauthor={#1}
}
\makeatother

\begin{document}

\begin{frontmatter}
\title{Monotone determined spaces via $\mathbb{C}$-generated spaces}
\tnotetext[t1]{Research supported by NSF of China (Nos. ).}

\author{Yuxu Chen}
\ead{chenyuxu@scu.edu.cn}
\author{Hui Kou}
\ead{kouhui@scu.edu.cn}
\cortext[cor]{Corresponding author}
\author{Zhenchao Lyu\corref{cor}}
\ead{zhenchaolyu@scu.edu.cn}

\address{Department of Mathematics, Sichuan University, Chengdu 610064, China}

\begin{abstract}
The category of monotone determined spaces is an extended topological framework for dcpos in domain theory. We first show that monotone determined spaces are exactly the spaces generated by one-point convergence spaces, and then naturally form a convenient Cartesian closed category of $\mathbb{C}$-generated spaces. 
We then show that monotone determined spaces are not always compact Hausdorff generated, answering the question raised by Ingo Battenfeld in 2013. Moreover, we generalize the notion of monotone determined spaces by introducing $\mathcal{C}$-determined spaces and showing that categories of $\mathcal{C}$-determined spaces correspond to coreflective subcategories of topological spaces. This yields a uniform construction of several convenient categories determined by directed, chain and monotone sequential convergence classes. We finally discuss the relationships among them, including the categories generated by continuous spaces, quasicontinuous spaces and Scott spaces of dcpos.\newline


\noindent {\bf Classification}: 18D15; 18F60; 54B30
\end{abstract}
\end{frontmatter}

\section{Introduction}

Monotone determined spaces were first introduced by Ern\'{e} in \cite{ERNE2009}. An interesting result for monotone determined spaces is that the lattice of all open sets is algebraic if and only if it is hyperalgebraic. Independently, by means of setting a final topology of a $T_{0}$ space, Yu and Kou \cite{YK2015} introduced directed spaces as a topological generalization of dcpos. In subsequent studies, it was found that directed spaces are equivalent to the monotone determined spaces. They also proved that $\bf Top_{md}$, the category of monotone determined spaces and continuous maps, is a coreflective Cartesian closed subcategory of the category of $T_0$ topological spaces, and contains all posets endowed with the Scott topology as a full subcategory. Luo and Xu independently proved the same results in \cite{XXQ2016}. In \cite{Battenfeld2013}, Ingo Battenfeld called these spaces crp-spaces and asked whether they are compact Hausdorff generated.

The category of monotone determined spaces is a natural topological extension \cite{YK2015,XXQ2016} of dcpos and provides a useful setting for domain-theoretic constructions. For example, free dcpo-algebras and free dcpo-cones can be obtained by taking the D-completion of free algebras over monotone determined spaces \cite{CKLX2024Cone}. Related constructions have been applied to upper powerdomains of quasicontinuous dcpos and free algebras over continuous spaces \cite{XK2020,CKL2024Upper,CKL2024Continuous}. Thus monotone determined spaces provide a route in which one first performs a construction topologically, such as a free monotone determined space algebra, and then returns to the classical dcpo setting by D-completion. In particular, for a dcpo $L$, the D-completion of the free algebra of $L$ in monotone determined spaces recovers the free algebra of $L$ in dcpos. These findings emphasize the importance of monotone determined spaces in domain theory and motivate the study of their properties.

Constructing $\mathbb{C}$-generated spaces \cite{HER72} is a convenient and powerful way to obtain coreflective Cartesian closed categories of topological spaces. In \cite{CG2004,LiMiao2023}, a series of $\mathbb{C}$-generated spaces related to domain theory, such as Kelley spaces, sequential spaces and qcb-spaces, are investigated. Monotone determined spaces, as natural topological extensions of dcpos, are more closely related to domain theory. Therefore, it is natural to ask whether monotone determined spaces can be described as $\mathbb{C}$-generated spaces. In this paper, we answer this question positively. We introduce a simple class of generating spaces, called one-point convergence spaces, and prove that the spaces generated by them are exactly the monotone determined spaces. We show that one-point convergence spaces are core compact and productive; therefore, monotone determined spaces form a convenient Cartesian closed category of $\mathbb{C}$-generated spaces. This provides a new perspective on monotone determined spaces and allows us to obtain many properties of them directly from the general theory of $\mathbb{C}$-generated spaces. Moreover, we show that monotone determined spaces are exactly locally hypercompact-space-generated. By contrast, spaces generated by continuous spaces are strictly contained in monotone determined spaces, and they can be described as spaces generated by strong one-point convergence spaces.

In \cite{CG2004}, it was shown that every Scott space of a dcpo is compact Hausdorff generated. As a natural extension of dcpos, it is reasonable to ask whether monotone determined spaces are compact Hausdorff generated. This question was also raised by Ingo Battenfeld when comparing different approaches for constructing free dcpo-algebras, and was left open in Section 6.2 of \cite{Battenfeld2013}. In this paper, we answer this question negatively and show that the compactly generated approach and the monotone determined approach are two different frameworks, neither of which contains the other. Compactly generated spaces and monotone determined spaces overlap in some important cases, such as dcpos, continuous spaces and quasicontinuous spaces, but they are different.

After these two results, we turn to the mechanism behind them. From the viewpoint of convergence classes, we introduce the notion of $\mathcal{C}$-determined spaces, where $\mathcal{C}$ is an operation assigning a convergence class to every topological space. We show that coreflective subcategories of topological spaces correspond to categories of $\mathcal{C}$-determined spaces. This gives a uniform way to construct categories determined by directed, chain and monotone sequential convergence classes; these categories are included in the monotone determined spaces and can also be described as Cartesian closed convenient categories of $\mathbb{C}$-generated spaces. We then compare the generated categories arising from continuous spaces, quasicontinuous spaces and Scott spaces of dcpos.

The paper is organized as follows. Section 2 recalls the basic facts on convergence classes and monotone determined spaces. Section 3 introduces one-point convergence spaces and proves that they generate exactly $\bf Top_{md}$. Section 4 answers Battenfeld's compact-generation question negatively and records the incomparability between $\bf Top_{md}$ and $\bf Top_{\mathbb K}$. Section 5 develops the general theory of $\mathcal{C}$-determined spaces and relates it to coreflective subcategories and $\mathbb{C}$-generated spaces. Section 6 compares the resulting generated categories, including the directed, chain and monotone sequential cases, the category generated by continuous spaces, the one generated by quasicontinuous spaces, and the category generated by Scott spaces of dcpos.

\section{Preliminaries}

We assume some basic knowledge of domain theory,\ topology and category theory,\ as in,\ e.g.,\ \cite{AJ94,CATE2004,GHKLMS03,GAU,GTP,Jung90,ZK2011}.

For any poset $(P,\leq)$ and $A\subseteq P$,\ we set
\[
A^{\ua} = \{x \in P : \forall a \in A,\ a \leq x \},\qquad
A^{\da} = \{ x \in P : \forall a \in A,\ x \leq a \}.
\]
Thus $A^{\ua}$ and $A^{\da}$ denote the sets of common upper and lower bounds of $A$. We also use
\[
\ua A=\{x\in P:\exists a\in A,\ a\leq x\},\qquad
\da A=\{x\in P:\exists a\in A,\ x\leq a\}
\]
for the generated upper and lower sets; for a singleton we write $\ua x$ and $\da x$. The cut of $A$ is defined by $A^\delta = A^{\ua\da} $. 
A nonempty subset $D$ of $P$ is directed iff every nonempty finite subset of\ $D$\ has an upper bound in\ $D$. We use\ $\sigma(P)$,\ $\upsilon(P)$ and $A(P)$ to denote the Scott topology,\ the upper topology,\ and the Alexandroff topology on $P$ respectively.\

Given any topological space $X$,\ its topology is denoted by $\mathcal{O}(X)$.\ For any net $\xi = \{x_j\}_{j\in J}$,\ where $J$ is a directed set,\ we write it $\{x_{j}\}$ for short.\ Given any $x\in X$, $\{ x_j \}$ is called converging to $x$,\ denoted by  $\{ x_j \} \rightarrow x $,\ if $\{ x_j \}$ is eventually in every open neighborhood of $x$.\ 
The specialization order on $X$ is defined by $x\sqsubseteq y$ iff $ x\in \overline{\{y\}}$. When $X$ is $T_{0}$,\ its specialization order is a partial order.\ Throughout this paper, all spaces are assumed to be $T_0$, and $\bf Top$ denotes the category of $T_0$ spaces and continuous maps unless otherwise stated.

The standard theory of $\mathbb C$-generated spaces will be used only in this
$T_0$ setting. This causes no difficulty: if $X$ is $T_0$, then the
$\mathbb C$-generated topology on $X$ is finer than the original topology, and
hence is again $T_0$. Thus the usual coreflective and Cartesian closed
arguments for $\mathbb C$-generated spaces restrict to the category of $T_0$
spaces in the form needed below.

For any poset $(P, \leq)$ endowed with a topology $\tau$,\ the specialization order $ \sqsubseteq$ of the topology $\tau$ agrees with the order $\leq$ iff $\tau$ is coarser than the Alexandroff topology and finer than the upper topology.\ We do not distinguish $\leq$ and $\sqsubseteq$ if they coincide.\

\vskip 3mm

It is convenient to define a topology through convergent nets.\ Some topological properties,\ like compactness, can be characterized by the convergent nets.\ We first illustrate some basic relationships between convergent nets and topologies, and then introduce the notion of monotone determined spaces.

Let $X$ be a set and $\Phi(X) = \{\xi: \xi \text{ is a net on } X\}$ be the class of all nets on $X$. A convergence class $\mathcal{E}$ on $X$ is defined as a subclass of $\Phi(X) \times X$. We also denote by $\xi \to_{\mathcal{E}} x$ if $(\xi,x) \in \mathcal{E}$. For any convergence class $\mathcal{E}$ on $X$, there is an induced topology $\mathcal{E}(X)$ as follows:\
 $$U \in \mathcal{E}(X) \Longleftrightarrow \forall (\xi,x) \in \mathcal{E},  x \in U \text{ implies } \xi \text{ is eventually in } U. $$ 
 The topology $\mathcal{E}(X)$ is said to be determined by $\mathcal{E}$.\ We denote by $\mathcal{E}X = (X,\mathcal{E}(X))$. 

In this paper, we mainly focus on the spaces determined by monotone convergence classes, i.e., convergence classes whose nets are directed sets under the specialization order. Restricting the nets to directed subsets gives the notion of a monotone determined space. Given any $T_0$ space $X$, define $$\mathcal{D}_X = \{(D,x): D \text{ is a directed subset of } X,\  D \to x \text{ in } X\}.$$

\begin{definition} \rm \cite{YK2015}
Let $(X,\tau)$ be a $T_0$ space. The topology determined by $\mathcal{D}_X$ is called the directed topology, denoted by $d(\tau)$. A space is called a directed space if $d(\tau) = \tau$.
\end{definition}

\begin{remark} \rm
In \cite{ERNE2009}, a topological space $X$ is said to be a monotone determined space if its topology is equal to the topology determined by its monotone convergent nets and limits.  The notion of directed spaces is in fact equivalent to monotone determined spaces. In \cite{YK2015}, Yu and Kou showed that the category of directed spaces is a coreflective Cartesian closed subcategory of $T_0$ topological spaces. Later, similar results were independently proved by Luo and Xu in \cite{XXQ2016}. In \cite{Battenfeld2013}, Ingo Battenfeld called monotone determined spaces crp-spaces. In the following, we uniformly use the name of monotone determined spaces.
\end{remark}

All monotone determined spaces together with continuous maps form a category,\ denoted by $\bf Top_{md}$. Given any space $(X,\tau)$, denote by $\mathcal{D}X = (X, d(\tau))$. The following are some basic properties of monotone determined spaces.

\begin{proposition} \rm  \cite{YK2015} \label{DIRE}
Let $(X,\tau )$ be a $T_0$ topological space. 
\begin{enumerate}
\item For any $U\in d(\tau)$,\ $U=\ua U$.
\item $X$ equipped with $d(\tau)$ is a $T_0$ topological space such that $\sqsubseteq_d =\sqsubseteq$,\ where $\sqsubseteq_d$ is the specialization order relative to $d(\tau)$.
\item For a directed subset $D$ of $X$,\ $D \rightarrow x$ iff $D\rightarrow_d x$ for all $x\in X$,\ where $D\rightarrow_d x$ means that $D$ converges to $x$ with respect to the topology $d(\tau)$.
\item  $d(d(\tau)) = d(\tau)$,\ i.e.,\ $\mathcal{D}X$ is a monotone determined space.\
\end{enumerate}
\end{proposition}

\begin{proposition} \rm \cite{XXQ2016,YK2015} \label{SEP}
Let $X,Y,Z$ be monotone determined spaces, and $X \otimes_d Y$ be the categorical products of $X$ and $Y$ in $\bf Top_{md}$. A map $f:X\otimes_d Y\longrightarrow Z$ is continuous iff it is separately continuous,\ i.e.,\ continuous with respect to each variable.
\end{proposition}

$\mathcal{D}$ can be viewed as the coreflector from $\bf Top_0$ to the coreflective subcategory $\bf Top_{md}$ as follows:
 \begin{eqnarray*}
& &\forall X\in {\rm Ob}({\bf Top}), \ \mathcal{D}(X) = \mathcal{D}X;\\
& &\forall f\in {\rm Mor}({\bf Top}), \ \mathcal{D}(f) (x) = f(x).
\end{eqnarray*}

Monotone determined spaces contain many interesting spaces in domain theory such as posets with the Scott topology,\ posets with the Alexandroff topology, continuous spaces and locally hypercompact spaces (see \cite{ERNE2009,XXQ2016,YK2015,HW2009,Li2009,Z2015}). 
\vskip 3mm

It is also easy to check the following lemma.

\begin{lemma} \label{con sup}
Let $X$ be a $T_{0}$ topological space and $D$ be a directed subset of $X$. 
\begin{enumerate}
\item[(1)] If $D \to x$, then $x \in D^\delta$.
\item[(2)] If $D \to x$ and $\sup D$ exists,\ then $x \leq \sup D$.\ 
\item[(3)] If $D\subseteq\da x$ and $D\rightarrow x$ for $x\in X$, then $x=\sup D$.
\end{enumerate}
\end{lemma}

For a topological space $X$ determined by a class $\mathcal{E}$, the continuous maps from $X$ can be characterized by the preservation for $\mathcal{E}$.

\begin{lemma} \label{CONV}
Let $Y$ be any topological space, $X$ be a set and $\mathcal{E}$ be a convergence class on $X$.
 A map $f$ from $(X,\mathcal{E}(X))$ into $Y$  is continuous iff for any $\xi \to_\mathcal{E} x$,\ $f(\xi) \rightarrow f(x)$ in $Y$ holds,\ i.e.,\ $f$ preserves  $\xi \rightarrow x$ for each $(\xi,x)$ in $\mathcal{E}$.
\end{lemma}

\begin{proof}
The necessity is obvious.\ For sufficiency,\ letting $U$ be an open subset of $Y$,\ we need only to show that $f^{-1}(U)$ is open in $(X,\mathcal{E}(X))$.\ For any $ (\xi,x) \in \mathcal{E}$, and $x \in f^{-1}(U)$,\ we have that $\xi$ is eventually in $f^{-1}(U)$.\ Otherwise,\ $f(\xi)$ is not eventually in $U$,\ a contradiction.
\end{proof}

We shall use the following standard classes of spaces and dcpos. A $T_0$ space
$X$ is called a $c$-space if, for every
$x\in U\in\mathcal{O}(X)$, there exists $y\in U$ such that
\(
x\in {\rm int}(\ua y),
\)
where $\ua y$ is taken with respect to the specialization order of $X$.
A space $X$ is called \textit{locally hypercompact} if, for
every $x\in U\in\mathcal{O}(X)$, there is a finite subset $F\subseteq X$ such
that
\(
x\in {\rm int}(\ua F)\subseteq \ua F\subseteq U.
\)
The notions of $c$-space and locally hypercompact space are also called
continuous space and quasicontinuous space, respectively~\cite{XK2022,CKL2024Continuous}. For a dcpo $P$, the
Scott space $\Sigma P=(P,\sigma(P))$ is a $c$-space iff $P$ is a continuous
dcpo; similarly, $\Sigma P$ is locally hypercompact iff $P$ is a
quasicontinuous dcpo \cite{GHKLMS03,GAU}. A space $X$ is called \textit{algebraic} if, for every $x\in U\in\mathcal{O}(X)$, there exists $k \in X$ such that
\(
x\in {\rm int}(\ua k) =  \ua k\subseteq U.
\)  A space $X$ is called \textit{quasialgebraic} if, for every $x\in U\in\mathcal{O}(X)$, there is a finite subset $F$ of $X$ such that
\(
x\in {\rm int}(\ua F) = \ua F\subseteq U.
\) Every algebraic space is a $c$-space, and every quasialgebraic space is a locally hypercompact space. 

A space $X$ is called \textit{locally compact} if, for every $x\in U\in\mathcal{O}(X)$,
there are an open set $V$ and a compact set $K$ such that
\(
x\in V\subseteq K\subseteq U.
\)
For open sets $U,V\in\mathcal{O}(X)$, write $U\ll V$ if every open cover of $V$
has a finite subfamily covering $U$. The space $X$ is called \textit{core compact} if,
for every $x\in V\in\mathcal{O}(X)$, there exists $U\in\mathcal{O}(X)$ such
that
\(
x\in U\ll V.
\)
Equivalently, $\mathcal{O}(X)$ is a continuous lattice \cite{GAU,LyuChenJia2022}. Every
locally hypercompact space is locally compact, and every locally compact space
is core compact.

\section{Monotone determined spaces as $\mathbb{C}$-generated spaces}

$\mathbb{C}$-generated spaces are special spaces that generate convenient Cartesian closed categories of topological spaces \cite{HER72}. Well known examples include Kelley spaces and sequential spaces.\ In \cite{CG2004}, a series of $\mathbb{C}$-generated spaces closely related to domain theory are investigated.  In this section, we introduce a new class of generating spaces, called one-point convergence spaces, and show that the collection of one-point convergence spaces generates exactly the monotone determined spaces. In particular, the collection of one-point convergence spaces is productive. Thus, monotone determined spaces form a convenient Cartesian closed category of $\mathbb{C}$-generated spaces. Many properties of monotone determined spaces can be obtained directly from this fact. 

\vskip 3mm

Let $\mathbb C$ be a collection of topological spaces,\ referred to as generating spaces.\ A map $k : C \to X$ is called a probe over $X$ if $C$ is in $\mathbb{C}$ and $k$ is continuous.\ A subset $U$ of a topological space $ X$ is called ${\mathbb C}$-open if its inverse image $k^{-1} (U) $ by any probe $k$ is open.\ All ${\mathbb C}$-open subsets of $X$ form a new finer topology on $X$,\ called the ${\mathbb C}$-generated topology\ (also called the final topology)\ on $X$ \cite{CG2004}.\ $X$ equipped with its ${\mathbb C} $-generated topology is denoted by ${\mathbb C} X$.\ A topological space $X$ is called a ${\mathbb C}$-generated space if $X={\mathbb C}X$. All ${\mathbb C}$-generated spaces together with continuous maps as morphisms form a category,\ denoted by $\bf Top_{\mathbb{C}}$.\
The following are some basic properties of coreflective subcategories of $\bf Top$.\

\begin{lemma}\rm  \cite{CG2004} \label{COLI}
A space is $\mathbb{C}$-generated iff it is a colimit in $\bf Top$ of generating spaces iff it is a quotient of disjoint sums of generating spaces.
\end{lemma}

A collection $\mathbb{C}$ of generating spaces is called productive if all topological spaces in $\mathbb{C}$ are core compact and every topological product $C_{1} \times C_{2}$,\ where $C_{1},C_{2} \in \mathbb{C}$, is $\mathbb{C}$-generated.

\begin{theorem} \rm  \cite{CG2004} \label{PROD}
If ${\mathbb C}$ is a productive class of topological spaces and contains at least one non-empty space $C_0$,\ then $\bf Top_{\mathbb C}
$ is Cartesian closed.
\end{theorem}

Now, we define one-point convergence spaces. For any directed poset $D$, let $\infty$ be an additional element, and let $\beta_D$ be the topology on $D \cup \{\infty\}$ generated by taking all $\ua x$ and $\ua x \cup \{\infty\}$ for $x \in D$ as a subbase.

\begin{definition}
Let $D$ be a directed poset. Then,  $(D \cup \{\infty\}, \beta_D)$ is called a one-point convergence space. Denote by $\mathbb{D}$ the collection of all one-point convergence spaces. 
\end{definition}

\begin{proposition} \label{RELA}
Given any topological space $X$,\ we have $\mathcal{D} X = \mathbb{D} X$.
\end{proposition}

\begin{proof}

(1) For any $(D \cup \{\infty\},\beta_{D})$ in $\mathbb{D}$,\ it is easy to check that $D \rightarrow \infty$ and $\{d\} \to e $ for any $e \leq d$ in $D$.\ Let $\mathcal{E}$ be the class of these pairs of convergent nets and limits,\ then it determines a topology $\mathcal{E}(D \cup \{\infty\})$ and we have $ \mathcal{E}(D \cup \{\infty\}) = \beta_{D}$.\ 
Given any open set $U \in \mathcal{E}(D \cup \{\infty\})$, and $d \in U$, since $ \{e\} \ra d$ for any $d \leq e$,\ we know $\ua d \subseteq U$.\ If $\infty \in U$,\ since $D \rightarrow \infty$,\ there exists some $d \in D$ such that $d \in U$,\ and then $\ua d \cup \{\infty\} \subseteq U$.\ Therefore, $U$ is the union of sets with the form $\ua d$ or $\ua d \cup \{\infty\} $.\ Thus,\ $ \beta_{D} = \mathcal{E}(D \cup \{\infty\})$.\
By Lemma \ref{CONV},\ for any map $p:(D \cup \{\infty\},\beta_{D}) \rightarrow X$,\ $p$ is a probe iff it preserves these limits in $\mathcal{E}$ iff $p$ is monotone and $p(D) \rightarrow p(\infty)$.\ So,\ if $p$ is a probe,\ then $p(D)$ is a directed subset in $X$ and $p(D) \rightarrow p(\infty)$.

On the other hand,\ given any directed subset $E$ of $X$ such that $E \rightarrow x$ in $X$,\ let $D=E$ with the induced specialization order, and let
$q:(D \cup \{\infty\},\beta_{D}) \to X$ be the map such that $q(d) = d$ for all $d \in D$ and $q(\infty) = x$.\ Then $q$ is continuous.\ Thus,\ $E \to x$ in the final topology $F(X)$.\

By definition,\ the final topology $F(X)$ on $X$ is the finest topology such that all probes are continuous,\ which is equivalent to that $F(X)$ is the finest topology such that if $D$ is directed and $D \rightarrow x$ in $X$,\ then $ D$ is directed and $D \rightarrow x$ in $(X,F(X))$ as well.\ Therefore,\ $F(X)$ is equal to the topology on $X$ determined by $\mathcal{D}_{X} = \{(D,x):D \text{ is a directed subset in } X, D \rightarrow x \text{ in } X\} $.\ 
\end{proof}

Since for any space $X$, $\mathcal{D} X$ is  a monotone determined space,\ and every monotone determined space $Y$ is equal to $\mathcal{D}Y$. We have the following conclusion. 

\begin{corollary}\label{d-generated-md}
Let $X$ be a topological space. Then, $X$ is $\mathbb{D}$-generated if and only if $X$ is monotone determined.
\end{corollary}

Therefore, one-point convergence spaces generate exactly the monotone determined spaces. Moreover, one-point convergence spaces are locally hypercompact, hence core compact. This is what is needed for productivity below. 

\begin{proposition} \rm  \label{monotone core}
One-point convergence spaces are quasialgebraic spaces. Consequently, they are core compact.
\end{proposition}
\begin{proof}
Given any one-point convergence space $X=(D \cup \{\infty\},\beta_{D})$,\ since $D$ is a directed set,\ it is easy to check that all $\ua x$ and $\ua x \cup \{\infty\}$ for $x \in D$ also form a base of $\beta_{D}$. We prove that $X$ is locally hypercompact.

Let $U \in \beta_{D}$ and $z \in U$. If $z \in D$, then $U$ is an upper set and hence $\ua z\subseteq U$. Taking $F=\{z\}$, we have
\(
z\in {\rm int}(\ua F)=\ua z\subseteq U.
\)

If $z=\infty$, then there is some $y\in D$ such that
\(
\ua y\cup\{\infty\}\subseteq U.
\)
Taking $F=\{y,\infty\}$, then
$\infty\in{\rm int}(\ua F) = \ua F\subseteq U$.
Thus $X$ is quasialgebraic.
\end{proof}

\begin{proposition}  \label{BRID}
 $\mathbb{D}$ is a productive class.
\end{proposition}

\begin{proof}
 
By Proposition \ref{monotone core}, every space in $\mathbb{D}$ is core compact.
Then, we need only to prove that every topological product of any two one-point convergence spaces is a monotone determined space.  

Let $X = (D_{1} \cup \{\infty_{1}\} , \beta_{D_{1}}) \times (D_{2} \cup \{ \infty_{2} \}, \beta_{D_{2}})$.\ By Proposition \ref{RELA},\ to show $X = \mathbb{D}X$,\ we need only to check that $X = \mathcal{D}X$.\ Obviously, $\{(a,b)\} \rightarrow (c,d)$ for every $(c,d) \sqsubseteq (a,b)$ in $X$,\ $\{(d,b)\}_{d \in D_{1}} \rightarrow (\infty_{1},b)$ for every $b \in D_{2} \cup \{\infty_{2}\}$,\ and $ \{(a,d)\}_{d \in D_{2}} \rightarrow (a,\infty_{2})$ for every $a \in D_{1} \cup \{\infty_{1}\} $.\ Conversely,\ we claim that these pairs of convergent nets and limits determines $X$. It suffices to show that the topology determined by these pairs, denoted by $\mathcal{E}(X)$, is coarser than $\mathcal{O}(X)$. Given any $U \in \mathcal{E}(X)$, it is an upper set since $\{(a,b)\} \rightarrow (c,d)$ for every $(c,d) \sqsubseteq (a,b)$ in $X$.\ Let $x=(a,b) \in U $.\ There are four cases:
   
  Case 1.\ $a \in D_{1}, b \in D_{2}$:\ $\ua (a,b) = \ua a \times \ua b \subseteq U$;\

  Case 2.\ $ a \in D_{1} , b = \infty_{2}$:\ by $\{(a,d)\}_{d \in D_{2}} \rightarrow (a,\infty_{2})$,\ we have that there exists a $d_{0} \in D_{2}$ such that $(a,d_{0}) \in U$. It follows that $\ua a \times (\ua d_{0} \cup \{\infty_{2}\}) \subseteq U$;\

  Case 3.\ $ a = \infty_{1} , b \in D_{2}$:\ similarly to case 2,\ there exists some $c_{0}$ such that $(\ua c_{0} \cup \{\infty_{1}\}) \times \ua b \subseteq U$;\

  Case 4.\ $x = (\infty_{1}, \infty_{2})$:\ by $\{(\infty_{1},d)\}_{d \in D_{2}} \rightarrow (\infty_{1},\infty_{2})$,\ there exists a $d_{0} \in D_{2}$ such that $(\infty_{1},d_{0}) \in U$.\ It follows that there exists a $c_{0} \in D_{1}$ such that $(c_{0},d_{0}) \in U$.\ Then $(\ua c_{0} \cup \{\infty_{1}\}) \times \ua d_{0} \subseteq U$.\ Similarly,\ there exist some $c_{1} \in D_{1}$,\ and  $d_{1} \in D_{2}$ such that $\ua c_{1} \times (\ua d_{1} \cup \{\infty_{2}\}) \subseteq U$.\ Let $c_{0},c_{1} \leq c_{2} \in D_{1}$  and $d_{0},d_{1} \leq d_{2} \in D_{2}$.\ We have $(\ua c_{2} \cup \{\infty_{1}\}) \times (\ua d_{2} \cup \{\infty_{2}\}) \subseteq U$.\

Thus, $\mathcal{E}(X) \subseteq \mathcal{O}(X)$.\ Since these pairs of convergent nets and limits of $\mathcal{E}$ are contained in $\mathcal{D}_{X}$,\ we have $\mathcal{O}(X) \subseteq \mathcal{D}(X) \subseteq \mathcal{E}(X) \subseteq \mathcal{O}(X)$, i.e, $X = \mathcal{D} X$.
 \end{proof}

 Then, by combining Theorem \ref{PROD} and Proposition \ref{BRID}, we know that  $\bf Top_{md}$ is a convenient Cartesian closed subcategory of topological spaces.

\begin{theorem}
$\mathbb{D}$-generated spaces are exactly monotone determined spaces, and $\bf Top_{md}$ is a convenient Cartesian closed subcategory of topological spaces.
\end{theorem}

By this observation, many properties of monotone determined spaces can be obtained directly. For example, $\bf Top_{md}$ is Cartesian closed and a coreflective subcategory of $T_0$ topological spaces.

\section{Monotone determined spaces are not compact Hausdorff generated}

Let $\bf Top_{\mathbb K}$ denote the category of compact Hausdorff generated spaces. Thus a space $X$ belongs to $\bf Top_{\mathbb K}$ iff a subset $U\subseteq X$ is open whenever $p^{-1}(U)$ is open for every continuous map $p:K\to X$ from a compact Hausdorff space $K$. In this section we answer negatively the question whether the category $\bf Top_{md}$ of monotone determined spaces is contained in $\bf Top_{\mathbb K}$. This settles, for monotone determined spaces, the compact-generation question left open in Section 6.2 of \cite{Battenfeld2013}.

To show that $\bf Top_{\mathbb D}\not\subseteq \bf Top_{\mathbb K}$, we need only to show that there exists a one-point convergence space in $\mathbb{D}$ that is not in $\bf Top_{\mathbb K}$. We now show that the one-point convergence space $X_D$ generated by the directed poset $D=[\mathbb R]^{<\omega}$ of all finite subsets of the real numbers is such a space. We first need the following lemma.

\begin{lemma}\label{CH-neighborhood-lemma}
Let $K$ be a compact Hausdorff space and let $b\in K$ be a non-isolated point, i.e., $\{b\}$ is not open. Suppose that $(N_i)_{i\in I}$ is a family of open neighbourhoods of $b$ such that, for every $x\neq b$, the set
\[
\{i\in I:x\in N_i\}
\]
is finite. Then $I$ is at most countable.
\end{lemma}

\begin{proof}
Suppose that $I$ is uncountable, and choose countably many distinct
$i_0,i_1,\ldots\in I$. Since every $x\neq b$ belongs to only finitely many of
the $N_{i_n}$, we have
\[
\bigcap_{n<\omega}N_{i_n}=\{b\}.
\]
Since $K$ is compact Hausdorff, it is regular and well-filtered. By regularity,
we can choose open sets $M_n$ such that
\[
b\in M_n\subseteq \overline{M_n}\subseteq N_{i_n}.
\]
Then the finite intersections
\[
B_n=M_0\cap\cdots\cap M_n
\]
form a countable local base at $b$. Indeed, let $O$ be an open neighbourhood of
$b$. The family
\[
K_n=\overline{B_n}\qquad(n<\omega)
\]
is a filtered family of compact saturated subsets of $K$, and
\[
\bigcap_{n<\omega}K_n
\subseteq
\bigcap_{n<\omega}N_{i_n}
=\{b\}\subseteq O.
\]
By well-filteredness, there exists $n<\omega$ such that $K_n\subseteq O$.
Hence $B_n\subseteq O$.

Now let \(\mathcal B=\{B_n:n<\omega\}\) be this countable local base. For every \(i\in I\), since \(N_i\) is a neighbourhood of \(b\), we can choose \(n(i)<\omega\) such that $B_{n(i)}\subseteq N_i.$

Fix \(n<\omega\). If infinitely many \(i\in I\) satisfy
\(
B_n\subseteq N_i,
\)
then, since \(b\) is non-isolated, we can choose
\(
x\in B_n\setminus\{b\}.
\)
Then \(x\) belongs to infinitely many \(N_i\), contradicting the hypothesis. Thus for each \(n\), the set
\[
\{i\in I:B_n\subseteq N_i\}
\]
is finite. Since
\[
I=\bigcup_{n<\omega}\{i\in I:B_n\subseteq N_i\},
\]
the set \(I\) is a countable union of finite sets. Hence \(I\) is countable,
contradicting the assumption that \(I\) is uncountable. Therefore \(I\) is at most countable.
\end{proof}

\begin{theorem}\label{directed-not-kgenerated}
Let $D=[\mathbb R]^{<\omega}$ denote the set of all finite subsets of the real numbers, ordered by inclusion, and let
\[
X_D=D\cup\{\infty\}
\]
carry the one-point convergence topology $\beta_D$, i.e., $\ua x$ and $\ua x \cup \{\infty\}$ for all $x \in D$ form a subbasis for the topology.
Then
\(
X_D\in \bf Top_{\mathbb D}\)
but \(
X_D\notin \bf Top_{\mathbb K}.
\) Consequently, \(
\bf Top_{\mathbb D}\not\subseteq \bf Top_{\mathbb K}.
\)
\end{theorem}

\begin{proof}
The poset $D=[\mathbb R]^{<\omega}$ is directed, so $X_D$ is in $\mathbb D$, $X_D\in \bf Top_{\mathbb D}$. We show that $\{\infty\}$ is $\mathbb{K}$-open but not open. It is not open because every neighbourhood of $\infty$ contains a basic neighbourhood of the form
\(
\ua_D S\cup\{\infty\}
\)
for some finite $S\subseteq\mathbb R$, and this set contains $S\in D$.

Let $p:K\to X_D$ be continuous, where $K$ is compact Hausdorff, and put $B=p^{-1}(\{\infty\})$. Since $D=\ua_D \emptyset$ is open in $X_D$, the set $K\setminus B=p^{-1}(D)$ is open, so $B$ is closed. For each $r\in\mathbb R$, set
\[
W_r=p^{-1}(\ua_D \{r\}\cup\{\infty\}).
\]
Then $W_r$ is open in $K$ and contains $B$. If $k\in K\setminus B$, then $p(k)$ is a finite subset of $\mathbb R$, and
\[
k\in W_r \Longleftrightarrow r\in p(k).
\]
Thus each $k\in K\setminus B$ belongs to only finitely many of the sets $W_r$.

Assume that $B$ is not open.  
Let
\[
q:K\to K/B
\]
be the quotient map, and let
\[
b=q(B).
\]
Since \(B\) is closed in the compact Hausdorff space \(K\), 
the quotient space obtained by collapsing \(B\) to a single point is again compact Hausdorff.
Because \(B\) is not open in \(K\), the point \(b\) is not isolated in \(K/B\).

For each \(r\in\mathbb R\), the set \(W_r\) is saturated with respect to \(q\),
because it contains the whole equivalence class \(B\). Hence
\[
N_r=q(W_r)
\]
is an open subset of \(K/B\). Moreover, \(N_r\) is a neighbourhood of \(b\).

Let \(y\in K/B\) with \(y\neq b\). Then \(y=q(k)\) for a unique
\(k\in K\setminus B\). We have
\[
y\in N_r
\iff
k\in W_r
\iff
r\in p(k).
\]
Since \(p(k)\) is finite, the set
\(
\{r\in\mathbb R:y\in N_r\}
\)
is finite.

Thus the family
\(
(N_r)_{r\in\mathbb R}
\)
is a family of open neighbourhoods of the non-isolated point \(b\) in the
compact Hausdorff space \(K/B\), and every point \(y\neq b\) belongs to only
finitely many of these neighbourhoods. This contradicts Lemma \ref{CH-neighborhood-lemma}, since \(\mathbb R\) is uncountable. Hence \(B\) must be open.

Therefore, $p^{-1}(\{\infty\})$ is open for every compact Hausdorff probe $p$, so $\{\infty\}$ is $K$-open. Since it is not open in $X_D$, the space $X_D$ is not compact Hausdorff generated.
\end{proof}

Conversely, it is easy to see that $\bf Top_{\mathbb K}$ is not contained in $\bf Top_{\mathbb D}$ either: every compact Hausdorff space is $\mathbb K$-generated, while a $T_{1}$ monotone determined space must be discrete, because its specialization order is equality and its monotone determined topology is therefore the discrete topology. Thus the compact Hausdorff generated category and the category of monotone determined spaces are incomparable.

Battenfeld~\cite{Battenfeld2013} asks whether the
compactly generated approach to domain theory can absorb the monotone determined approach. This result shows that they are two frameworks which are not contained in one another.

\begin{remark}
There also exist non-$T_1$ spaces that are compact Hausdorff generated but not monotone determined. For example, let
$\mathbb S=\{\bot,\top\}$ be the Sierpi\'nski space, whose open subsets are
$\emptyset,\{\top\}$ and $\mathbb S$, and put
\[
X=\mathbb S\times [0,1].
\]
Then $X$ is $T_0$ but not $T_1$. It is compact Hausdorff generated. Indeed, the
map
\[
q:[0,1]\longrightarrow \mathbb S,\qquad
q(t)=\top\ (t<1),\qquad q(1)=\bot,
\]
is a quotient map. Since products with locally compact Hausdorff spaces
preserve quotient maps, the map
\[
q\times {\rm id}_{[0,1]}:[0,1]^2\longrightarrow \mathbb S\times[0,1]
\]
is a quotient map. Hence $X$ is a quotient of the compact Hausdorff space
$[0,1]^2$, and therefore $X\in\bf Top_{\mathbb K}$.

However, $X$ is not monotone determined. The specialization order of $X$ is
given by
\[
(i,t)\leq (j,s)
\quad\Longleftrightarrow\quad
i\leq j\text{ in }\mathbb S\ \text{and}\ t=s.
\]
Thus every directed subset of $X$ lies in a single vertical fibre
$\mathbb S\times\{t\}$. Consider
\[
U=\{\top\}\times\{0\}.
\]
This set is not open in $X$, since $\{0\}$ is not open in $[0,1]$. We show that
it is open in the directed topology. Let $D$ be a directed subset of $X$ with
$D\to x\in U$. Then $x=(\top,0)$. Since $D$ lies in some fibre
$\mathbb S\times\{t\}$ and the second projection is continuous, convergence to
$(\top,0)$ forces $t=0$. Moreover, $\{\top\}\times[0,1]$ is an open
neighbourhood of $(\top,0)$, so $D$ is eventually in
$\{\top\}\times[0,1]$. Hence $D\cap U\neq\emptyset$, and since $U$ is an upper
set, $D$ is eventually in $U$. Therefore $U$ is monotone determined open but not open in
$X$, so $X\notin\bf Top_{\mathbb D}$.
\end{remark}

\section{$\mathcal{C}$-determined spaces and generated spaces}

\subsection{Convergence class operations}

The preceding sections show that monotone determined spaces fit naturally into the theory of $\mathbb{C}$-generated spaces. We now isolate the common convergence-class pattern behind this phenomenon. We define the notion of a convergence class operation $\mathcal{C}$, which maps any topological space $X$ to a convergence class $\mathcal{C}_X$, and then there is a natural topological space determined by $\mathcal{C}_X$. Such spaces are called $\mathcal{C}$-determined spaces. Monotone determined spaces are concrete examples of $\mathcal{C}$-determined spaces obtained by setting $\mathcal{C}_X = \mathcal{D}_X$. We give an equivalent characterization of coreflective subcategories of topological spaces from the viewpoint of convergence classes: a category is a coreflective subcategory of $\bf Top_0$ iff it is the category of $\mathcal{C}$-determined spaces for some idempotent and consistent $\mathcal{C}$.  In particular, we investigate several concrete $\mathcal{C}$-determined spaces contained in monotone determined spaces and show that, on the other hand, they can be described as different kinds of $\mathbb{C}$-generated spaces. By the tool of $\mathbb{C}$-generated spaces, we can prove their Cartesian closedness uniformly. 

The following theorem gives a characterization of coreflective subcategories of $\bf Top$.\ In our $T_0$ convention, quotients and colimits are taken in $\bf Top$, i.e. with $T_0$-reflection when necessary \cite{XuShenXiZhao2020,XuZhao2020Rudin,Xu2021HSober,MiaoWangLi2024,MiaoJiaShenLi2024,LiLiZhao2025}.\

\begin{theorem} \cite{HER72}
A subcategory of $\bf Top$ is coreflective in $\bf Top$ iff it is invariant under the formation of disjoint topological unions and topological quotient spaces in $\bf Top$.\

Every coreflective subcategory of $\bf Top$ is cocomplete and is a cocomplete subcategory of $\bf Top$.

Every coreflective subcategory $\bf C$ of $\bf Top$ is complete and the limit in $\bf C$ is the $\bf C$-coreflection for the limit in $\bf Top$.\
\end{theorem}

A convergence class is called topological if it is the class of all pairs of convergent nets and their limits of a topological space.
For a topological space $X$, we denote by $\mathcal{T}_{X}$ the class of all pairs of convergent nets and their limits in $X$.

\begin{lemma} \label{coarser}
Let $X$ be a set and $\mathcal{E}, \mathcal{F}$ be convergence classes on $X$. If  $\mathcal{E} \subseteq \mathcal{F}$, then $(X,\mathcal{F}(X))$ is coarser than $(X,\mathcal{E} (X))$. 
\end{lemma}
\begin{proof}
Let $U$ be an open subset in $(X,\mathcal{F}(X))$,\ then for any $(\xi,x) \in \mathcal{F}$ and $x \in U$,\ $\xi$ is eventually in $U$.\ If $\mathcal{E} \subseteq \mathcal{F}$,\ we have that for any $(\xi,x) \in \mathcal{E}$ and $x \in U$, $\xi$ is eventually in $U$,\ i.e,\ $U$ is open in $(X,\mathcal{E}(X))$.
\end{proof}

\begin{definition}
We call $\mathcal{C}$ a convergence class operation if, for every topological space $X$, it assigns a convergence class $\mathcal{C}_{X}$ with
\(
\mathcal{C}_{X}\subseteq \mathcal{T}_{X}.
\)
Thus every convergence pair selected by $\mathcal{C}_{X}$ is an actual convergent net and limit in $X$.\ For such an operation, denote $(X,\mathcal{C}_{X}(X))$ by $\mathcal{C}X$.\ If $\mathcal{C}_{X} \subseteq \mathcal{C}_{\mathcal{C}X}$ holds for every topological space $X$, then we say that $\mathcal{C}$ is idempotent.\ $X$ is called a $\mathcal{C}$-determined space iff $X = \mathcal{C} X$.
\end{definition}

\begin{lemma} \label{C-determined}
Let $X$ be a topological space. 
\begin{enumerate}

\item[(1)] For any idempotent operation $\mathcal{C}$,\ $\mathcal{C} X$ is a $\mathcal{C}$-determined space. 

\item[(2)] For any two idempotent operations $\mathcal{C}, \mathcal{G}$,\ if $\mathcal{C}_{Y} \subseteq \mathcal{G}_{Y} $ holds for every topological space $Y$ (in which case $\mathcal{G}$ is called larger than $\mathcal{C}$),\ then every $\mathcal{C}$-determined space is a ${\mathcal{G}}$-determined space. 

\end{enumerate}
\end{lemma}

\begin{proof}

(1) For any topological space $X$,\ since $\mathcal{C}_{X} \subseteq \mathcal{C}_{\mathcal{C}X}$,\ Lemma \ref{coarser} implies that $ \mathcal{C}\mathcal{C} X $ is coarser than $\mathcal{C}X$.\ Conversely, by the definition of a convergence class operation,
\(
\mathcal{C}_{\mathcal{C}X}\subseteq \mathcal{T}_{\mathcal{C}X}.
\)
Hence every open set of $\mathcal{C}X$ is open with respect to all pairs in $\mathcal{C}_{\mathcal{C}X}$, so $\mathcal{C}X$ is coarser than $\mathcal{C} \mathcal{C}X$.\ Therefore $\mathcal{C} X$ is a $\mathcal{C}$-determined space. 

(2) Assuming that $Y$ is a $\mathcal{C}$-determined space,\ then $ \mathcal{C}Y = Y$ and, since $\mathcal{G}_{Y}\subseteq\mathcal{T}_{Y}$, the space $Y$ is coarser than $\mathcal{G}Y$.\ By the assumption that $\mathcal{C}_{Y} \subseteq \mathcal{G}_{Y}$ and Lemma \ref{coarser},\ we have that $\mathcal{G}Y$ is coarser than $ \mathcal{C}Y = Y$.\ Thus,\ $Y$ is also a $\mathcal{G}$-determined space.
\end{proof}

\begin{definition} Let $X,Y$ be two topological spaces and $\mathcal{C}$ be an idempotent operation.\ A map $f: X \to Y$ is said to be $\mathcal{C}$-continuous if $f: \mathcal{C}X \to  \mathcal{C}Y$ is continuous.
\end{definition}

\begin{definition} A convergence class operation $\mathcal{C}$ is said to be consistent if for any two topological spaces $X,Y$ and any continuous map $f:X \to Y$,\ $(\xi,x) \in \mathcal{C}_{X}$ implies $(f(\xi),f(x)) \in \mathcal{C}_{Y}$.
\end{definition}

\begin{lemma} \label{PCON}
 Let $\mathcal{C}$ be an idempotent and consistent operation.\ Given any topological space $Y$ and any $\mathcal{C}$-determined space $X$,  a map $f: X \rightarrow Y$ is continuous iff $f:X\rightarrow  \mathcal{C}Y$ is $\mathcal{C}$-continuous.
\end{lemma}

\begin{proof}
Suppose that $f: X \to Y$ is continuous. Since $\mathcal{C}$ is consistent,
\[
(\xi,x)\in\mathcal C_X
\quad\Longrightarrow\quad
(f(\xi),f(x))\in\mathcal C_Y.
\]
By Lemma \ref{CONV}, $f:\mathcal C X\to\mathcal C Y$ is continuous. Hence
$f:X\to\mathcal C Y$ is continuous. Since $\mathcal C\mathcal C Y=\mathcal C Y$,
we have that $f:X\to\mathcal C Y$ is $\mathcal C$-continuous.

For the converse, suppose that $f:X\to\mathcal C Y$ is $\mathcal C$-continuous.
Then $f:\mathcal C X\to\mathcal C\mathcal C Y$ is continuous. Since
$\mathcal C X=X$ and $\mathcal C\mathcal C Y=\mathcal C Y$ by Lemma
\ref{C-determined}, the map $f:X\to\mathcal C Y$ is continuous. Since $Y$ is
coarser than $\mathcal C Y$, the map $f:X\to Y$ is continuous.\
\end{proof}

\begin{corollary}\label{PCO}
Let $\mathcal{C}$ be an idempotent and consistent operation and $X,Y$ be $\mathcal{C}$-determined spaces.\ A map $f:X \to Y$ is continuous iff it is $\mathcal{C}$-continuous.\
\end{corollary}

For any idempotent and consistent operation $\mathcal{C}$,\ we denote $\bf Top_{\mathcal{C}}$ the category of all $\mathcal{C}$-determined spaces with $\mathcal{C}$-continuous maps as morphisms.\ Since a map between $\mathcal{C}$-determined spaces is continuous iff it is $\mathcal{C}$-continuous, $\bf Top_{\mathcal{C}}$ is a full subcategory of $\bf Top$.\ 

We now define $\widetilde{\mathcal{C}}$ as follows:
 \begin{eqnarray*}
& &\forall X\in {\rm Ob}({\bf Top}), \ \widetilde{\mathcal{C}}(X) = \mathcal{C}X,\\
& &\forall f\in {\rm Mor}({\bf Top}), \ \widetilde{\mathcal{C}}(f) (x) = f(x).
\end{eqnarray*}

Because $\mathcal{C}X$ is a $\mathcal{C}$-determined space for any topological space $X$, and the continuity is equivalent to $\mathcal{C}$-continuity for maps between $\mathcal{C}$-determined spaces,\ we know that $\widetilde{\mathcal{C}}$ is a functor from $\bf Top$ to $\bf Top_{\mathcal{C}}$.\ Moreover,\ for every idempotent and consistent $\mathcal{C}$,\ we show that $\bf Top_{\mathcal{C}}$ is a coreflective subcategory of $\bf Top$ \ and $\widetilde{\mathcal{C}}$ is the coreflector.\

\begin{definition}\cite{CATE2004}
Let $\bf A$ be a subcategory of $\bf B$ and let $B$ be a $\bf B$-object.\
\begin{enumerate}
\item An $\bf A$-coreflection for $B$ is a $\bf B$-morphism $A \overset{c} \to B$ from a $\bf A$-object $A$ to $B$ with the following universal property:\ for any $\bf B$-morphism $A^{\prime} \overset{f} \to B$ from some $\bf A$-object $A^{\prime}$ to $B$ there exists a unique $\bf A$-morphism $f^{\prime} : A^{\prime} \to A$ such that $f = c \circ  f^{\prime}$.\ Such $A$ is called a $\bf A$-coreflection for $B$.\

\item $\bf A$ is called a coreflective subcategory of $\bf B$ provided that each $\bf B$-object has an $\bf A$-coreflection.
\end{enumerate}
\end{definition}

\begin{theorem}\label{SUBCA}
Let $\mathcal{C}$ be an idempotent and consistent operation. Then $\bf Top_{\mathcal{C}}$ is a coreflective subcategory of $\bf Top$.\
\end{theorem}

\begin{proof}
Given any $Y$ in $\bf Top$,\ the identity map $id : \mathcal{C}Y \to Y$ is continuous.\ For any $X$ in $\bf Top_{\mathcal{C}}$ and continuous map $f:X \to Y$,\ $f: X \to \mathcal{C}Y$ is continuous by Lemma \ref{PCON} and Corollary \ref{PCO}.\ Therefore,\ $\mathcal{C} Y$ is the $\bf Top_{\mathcal{C}}$-coreflection for $Y$.\ 
\end{proof}
 
\begin{example} \label{EP}
These are some examples of $\mathcal{C}$-determined spaces for idempotent and consistent $\mathcal{C}$.
\begin{enumerate}
\item Discrete spaces are $1$-determined, where $1_{X} = \{(\{x\},x):x \in X\}$ for any space $X$. 
\item Alexandroff spaces are $ \mathcal{S}$-determined, where $\mathcal{S}_{X} = \{(\{y\},x):x,y \in X, x \sqsubseteq y\}$. 
\item Sequential spaces are $\mathcal{N}^{\prime \prime}$-determined, where
\[
\mathcal{N}^{\prime \prime}_{X}
=\{(M,x):x\in X,\ M \text{ is a sequence in } X,\ M \to x \}.
\]

\end{enumerate}
\end{example}

Conversely,\ given any coreflective subcategory $\bf C$ of $\bf Top$,\ we can find an idempotent and consistent operation $\mathcal{C}$ such that $\bf C = \bf Top_{\mathcal{C}}$.\ 

\begin{proposition} \cite{CATE2004} \label{CO}
If $\bf A$ is a coreflective subcategory of $\bf B$ and for each $\bf B$-object $B$,\ $A_{B} \overset{c_{B}} \rightarrow B$ is an $\bf A$-coreflection,\ then there exists a unique functor $C: \bf B \to \bf A$ (called a coreflector for $\bf A$) such that the following conditions are satisfied:
\begin{enumerate}
\item $C(B) = A_{B}$ for each $\bf B$-object $B$,\
\item for each $\bf B$-morphism $f:B \to B^{\prime}$ the diagram below commutes.\

$$\xymatrix{
  C(B) \ar[r]^{c_{B}} \ar[d]_{C(f)} &  B \ar[d]^{f}    \\
  C(B^{\prime}) \ar[r]_{c_{B^{\prime}}} &  B^{\prime} }$$
\end{enumerate}
\end{proposition}

\begin{theorem} \cite{HER72} \label{COREF}
A subcategory $\bf C$ of $\bf Top$ is coreflective in $\bf Top$ iff for each space $(X,\tau)$ in $\bf Top$,\ there exists a topology $\eta_{X}$ on $X$ with the following properties:
\begin{enumerate}
\item $\eta_{X}$ is finer than $\tau$
\item $(X,\eta_{X})$ is in $\bf C$
\item $\eta_{X}$ is the coarsest topology on $X$ satisfying $(1)$ and $(2)$
\item for each space $(Y,\phi)$ and each continuous map $f:(X,\tau) \to (Y,\phi)$,\ the same set map $f:(X,\eta_{X}) \to (Y,\eta_{Y})$ is continuous.\
\end{enumerate}
\end{theorem}

From Proposition \ref{CO} and Theorem \ref{COREF},\ we know that $C(X)$,\ the $\bf C$-coreflection for $(X,\tau)$ in $\bf Top$,\ is $(X,\eta_{X})$ up to isomorphism.\ Then we have the following statement.\

\begin{proposition}\label{REV}
Let $\bf C$ be a coreflective subcategory of $\bf Top$.\
Taking $\mathcal{C}$ an operation such that $\mathcal{C}_{X} = \mathcal{T}_{C(X)}$, where $\mathcal{T}_{C(X)}$ is the topological convergence class of $C(X)$,\ then $\mathcal{C}$ is idempotent and consistent.\ Moreover,\ for any topological space $X$,\ $C(X) = \mathcal{C}X$,\ i.e.,\ $\widetilde{\mathcal{C}}$ is the coreflector for $\bf C$.\
\end{proposition}

\begin{proof} Since the coreflection map $C(X)\to X$ is continuous, $\mathcal{T}_{C(X)}\subseteq\mathcal{T}_X$; hence $\mathcal{C}$ is a convergence class operation.\ By Theorem \ref{COREF},\ given any two topological spaces $X_{1},X_{2}$,\ if $X_{1}$ is coarser than $X_{2}$,\ then $ \mathcal{C}_{X_{2}} =\mathcal{T}_{C(X_{2})}  \subseteq   \mathcal{T}_{C(X_{1})} = \mathcal{C}_{X_{1}}$.\ Since $\mathcal{C}_{X} = \mathcal{T}_{C(X)}$,\ we have $\mathcal{C}X = (X,\mathcal{C}_{X}(X)) = (X,\mathcal{T}_{C(X)}(X)) = C(X)$.\ Obviously,\ for any space $X$ in $\bf C$,\ $C(X) = X$ and then $\mathcal{C}X = C(X) = X$.\ Thus,\ $\mathcal{C}_{X} = \mathcal{C}_{\mathcal{C}X} $,\ that is,\ $\mathcal{C}$ is idempotent.\ Given any continuous map $f:X \to Y$ in $\bf Top$,\ by Theorem \ref{COREF},\ $f: C(X) \to C(Y)$ is continuous.\ Then $(\xi,x) \in \mathcal{T}_{C(X)}$ implies $(f(\xi),f(x)) \in \mathcal{T}_{C(Y)}$,\ i.e.,\ $\mathcal{C}$ is consistent.\ 
\end{proof}

By Proposition \ref{SUBCA} and Proposition \ref{REV},\ we know that there is a correspondence between idempotent and consistent operations and coreflective subcategories.\ We define an equivalence relation on idempotent and consistent operations by declaring two operations $\mathcal{C},\mathcal{G}$ equivalent iff $\mathcal{G}X =\mathcal{C}X$ for every topological space $X$.\ The equivalence class of $\mathcal{C}$ is denoted by $[\mathcal{C}]$.\ It is easily seen that there is a largest operation among all operations in $[\mathcal{C}]$,\ namely the operation that assigns $\mathcal T_{\mathcal C X}$ to each space $X$.\ This is still contained in $\mathcal T_X$, because the coreflection map $\mathcal C X\to X$ is continuous.\
Then there is a one-to-one correspondence between all equivalent classes of idempotent and consistent operations and coreflective subcategories.\

It is easily seen that for any two idempotent and consistent operations $\mathcal{C},\mathcal{G}$,\ if $\mathcal{G}$ is larger than $\mathcal{C}$,\ i.e.,\ $ \mathcal{C}_{X} \subseteq \mathcal{G}_{X}$ for every topological space $X$,\ then $\bf Top_{\mathcal{C}}$ is a coreflective subcategory of $\bf Top_{\mathcal{G}}$.\

\subsection{Concrete operations and one-point generators}

Next\ we introduce some other concrete examples of $\mathcal{C}$-determined spaces as follows, all of which are contained in monotone determined spaces, and form convenient Cartesian closed categories via $\mathbb{C}$-generated spaces. The first case is just the monotone determined spaces; we include it here for uniformity.

\begin{definition}
We define six operations $\mathcal{D},\mathcal{D}^{\prime},\mathcal{I},\mathcal{I}^{\prime},\mathcal{N},\mathcal{N}^{\prime}$ as follows.\ Let $X$ be any topological space.\
\begin{enumerate}
\item $ \mathcal{D}_{X} = \{(D,x): D \text{ is directed in } X,\ D \rightarrow x \} $;
\item $\mathcal{D}^{\prime}_{X}= \{(D,x): D \text{ is directed in } X,\ x \in \da D \text{ or }(D \rightarrow x \text{ and } D \leq x)\} $;
\item $\mathcal{I}_{X} = \{(I,x): I \text{ is a chain in the specialization order of } X,\ I \rightarrow x \}$;
\item $\mathcal{I}^{\prime}_{X} = \{(I,x): I \text{ is a chain in the specialization order of } X,\ x \in \da I \text{ or }(I \rightarrow x \text{ and } I \leq x) \}$;
\item $\mathcal{N}_{X}= \{(M,x): M \text{ is a monotone sequence in } X,\ M \rightarrow x \}$;
\item $\mathcal{N}^{\prime}_{X}= \{(M,x): M \text{ is a monotone sequence in } X,\ x \in \da M \text{ or }(M \rightarrow x \text{ and } M \leq x) \}$.
\end{enumerate}
\end{definition}

The six operations above are all idempotent and consistent.\ We use operations $\mathcal{C}_{i}\ (1\leq i \leq 6)$ to denote these operations for convenience.\ By definition,\ $\mathcal{D}$-determined spaces are just the monotone determined spaces.\ Spaces determined by operations $\mathcal{C}_{i}\ (2\leq i \leq 6)$ are all included in $\mathcal{D}$-determined spaces.\ They all form coreflective subcategories of $\bf Top$ and then are cocomplete and complete.\ If we denote $\omega\mathcal{D}_{X}= \{(D,x) \in \mathcal{D}_{X} : \text{ D is countable }\}$, then it is easy to see that given any topological space $X$, $\omega\mathcal{D}X = \mathcal{N}X$.

Similarly to one-point convergence spaces,\ we define classes of generating spaces such that the spaces generated by them coincide respectively with the spaces determined by $\mathcal{C}_{i}\ (1\leq i\leq 6)$.

\begin{definition} \label{GENERATE}
Let $\mathscr{D}$ denote the class of all directed posets,\ $\mathscr{I}$ denote the class of all chains and $N$ denote the poset of natural numbers.\ Let $\mathscr{D}^{\ast}$ and $\mathscr{I}^{\ast}$ denote respectively the subclasses of directed posets and chains with no largest element.\ Given any poset $P$,\ let $\infty$ be an element that is not in $P$.\ We define two topologies on $P \cup \{\infty\}$ as follows:\ $\beta_{P}$ is the topology by taking all $\ua x $ and $ \ua x \cup \{\infty\} $ for $x \in P$ as a subbase;\ $\gamma_{P}$ is the topology by taking all $\ua x \cup \{\infty\}$ for $x \in P$ as a subbase.\ When $P$ is directed and has no largest element, the spaces $(P\cup\{\infty\},\gamma_P)$ are called strong one-point convergence spaces. The no-largest restriction keeps these generators $T_0$ and does not change the generated topologies, since any directed or chain convergence can be reindexed over a no-largest directed poset or chain. In particular,\ we define six classes of generating spaces as follows:
\begin{enumerate}

\item[(1)] $ \mathbb{D} = \{(D \cup \{\infty\},\beta_{D}): D \in  \mathscr{D} \} $;
\item[(2)] $\mathbb{D}^{\prime}= \{(D \cup \{\infty\},\gamma_{D}): D \in  \mathscr{D}^{\ast}\} $;
\item[(3)] $\mathbb{I} = \{(I \cup \{\infty\},\beta_{I}): I \in  \mathscr{I} \}$;
\item[(4)] $\mathbb{I}^{\prime} = \{(I \cup \{\infty\},\gamma_{I}): I \in  \mathscr{I}^{\ast} \}$;
\item[(5)] $\mathbb{N}= \{(N \cup \{\infty\},\beta_{N}) \}$;
\item[(6)] $\mathbb{N}^{\prime} = \{(N \cup \{\infty\},\gamma_{N}) \}$.

\end{enumerate}
\end{definition}

For convenience, we denote by $\mathbb{C}_{i}\ (1 \leq i \leq 6)$ the six classes in Definition \ref{GENERATE}. Then we have the following statement. The proof is analogous to the proof of Proposition \ref{RELA}, so we omit it.

\begin{proposition} \label{RELA1}

Given any topological space $X$,\ we have $\mathcal{C}_{i} X = \mathbb{C}_{i} X$ for $1 \leq i\leq 6$.

\end{proposition}

\begin{proposition}\label{Ci-productive}
For each $1\leq i\leq 6$, the class $\mathbb{C}_i$ is productive.
\end{proposition}

\begin{proof}
The generators in the $\beta$-families are core compact by Proposition
\ref{monotone core}. The same argument applies to the $\gamma$-families: if
$P$ is a no-largest directed poset and $X=(P\cup\{\infty\},\gamma_P)$, then every basic
neighbourhood of a point contains a principal saturated neighbourhood
$\ua p\cup\{\infty\}$ for some $p\in P$, so $X$ is core compact.

It remains to check products of generators. For the $\beta$-families, let
$\mathscr P$ denote one of $\mathscr D$, $\mathscr I$ and $\{N\}$; for the
$\gamma$-families, let $\mathscr P$ denote one of $\mathscr D^\ast$,
$\mathscr I^\ast$ and $\{N\}$. Let $\theta$ be the corresponding choice
$\beta$ or $\gamma$. Take
\[
X_P=(P\cup\{\infty_P\},\theta_P),\qquad
X_Q=(Q\cup\{\infty_Q\},\theta_Q),
\]
where $P,Q\in\mathscr P$. In the product $X_P\times X_Q$, consider the
singleton specialization convergences and the coordinate convergences
\[
\{(p,v)\}_{p\in P}\to(\infty_P,v)\quad (v\in Q\cup\{\infty_Q\}),
\]
and
\[
\{(u,q)\}_{q\in Q}\to(u,\infty_Q)\quad (u\in P\cup\{\infty_P\}).
\]
For $\theta=\gamma$ these coordinate convergences are also bounded above by
their limits, so they belong to the corresponding primed operation
$\mathcal C_i$; for $\theta=\beta$ they belong to the corresponding unprimed
operation.

These convergence pairs determine the product topology. Indeed, any subset
which is open for these pairs is an upper set. If it contains a point
$(p,q)$ with $p\in P$ and $q\in Q$, it contains the upper rectangle
$\ua p\times\ua q$. If it contains $(p,\infty_Q)$, the second coordinate
convergence gives some $q_0\in Q$ in the set, and upperness gives a basic
rectangle around $(p,\infty_Q)$ contained in the set; the case
$(\infty_P,q)$ is symmetric. Finally, if it contains
$(\infty_P,\infty_Q)$, applying the two coordinate convergences successively
and then using directedness of $P$ and $Q$ gives $p_0\in P$ and $q_0\in Q$
such that the basic rectangle
\[
(\ua p_0\cup\{\infty_P\})\times(\ua q_0\cup\{\infty_Q\})
\]
is contained in the set. Thus every set open for these convergence pairs is
product-open.

Hence $X_P\times X_Q=\mathcal C_i(X_P\times X_Q)$ for the corresponding
operation $\mathcal C_i$, and Proposition \ref{RELA1} gives
$X_P\times X_Q=\mathbb C_i(X_P\times X_Q)$. Therefore every binary product of
generators is $\mathbb C_i$-generated, and $\mathbb C_i$ is productive.
\end{proof}

\begin{theorem}\label{Ci-cc}
For each $1\leq i\leq 6$, the category $\bf Top_{\mathbb C_i}$ is Cartesian closed.
\end{theorem}

\begin{proof}
This follows immediately from Theorem \ref{PROD} and Proposition
\ref{Ci-productive}.
\end{proof}

\subsection{Continuous, quasicontinuous and dcpo Scott generators}

We now compare the one-point generators above with more familiar classes from
domain theory. The aim of this subsection is to identify the generated
categories obtained from continuous spaces, from quasicontinuous
spaces, and from Scott spaces of dcpos. The results below show that continuous
spaces generate the smaller category $\bf Top_{\mathbb D'}$, while
quasicontinuous spaces generate all monotone determined spaces. By
contrast, Scott spaces of algebraic, continuous, quasicontinuous and arbitrary
dcpos all generate the same category, which sits properly inside
$\bf Top_{\mathbb D'}$.

\begin{proposition}\label{gamma-cspaces}
For every directed poset $D$ with no largest element, the strong one-point convergence space $(D\cup\{\infty\},\gamma_D)$ is an algebraic space, hence also a continuous space.
\end{proposition}

\begin{proof}
The specialization order of $(D\cup\{\infty\},\gamma_D)$ is the order on $D$ with $\infty$ added as a largest element. Moreover, for every $d\in D$, \( \ua d=\ua_D d\cup\{\infty\} \) is open. Let $U$ be an open neighbourhood of a point $x$. If $x\in D$, then, since $U$ is an upper set and $\ua x$ is open, we have \( x\in{\rm int}(\ua x)=\ua x\subseteq U. \) If $x=\infty$, and $U$ is an open neighbourhood of $\infty$, then there exists some $d \in D$ such that
\(\infty\in{\rm int}(\ua d)=\ua d\subseteq U. 
\) Hence it is an algebraic space.
\end{proof}

Denote by $\mathsf{ConS}$ the class of all continuous spaces and by $\mathsf{QConS}$ the class of all quasicontinuous spaces. Denote by $\mathsf{AlgS}$ the class of all algebraic spaces and by $\mathsf{QAlgS}$ the class of all quasialgebraic spaces.

The following theorem shows that the category of continuous space-generated spaces coincides with the category of strong one-point convergence space-generated spaces.

\begin{theorem}\label{cspaces-dprime-generated}

\(
\bf Top_{\mathsf{AlgS}} = \bf Top_{\mathsf{ConS}}=\bf Top_{\mathbb D'}.
\)
\end{theorem}

\begin{proof}
By Proposition \ref{gamma-cspaces}, every generator in $\mathbb D'$ is an algebraic space. 
Hence $\bf Top_{\mathbb D'} \subseteq \bf Top_{\mathsf{AlgS}} \subseteq \bf Top_{\mathsf{ConS}}$.

Conversely, we show that every continuous space belongs to
$\bf Top_{\mathbb D'}$. Let $X$ be a continuous space and let
$U\in\mathcal D'(X)$. We prove that $U$ is open in the original topology of
$X$.

First $U$ is an upper set. Indeed, if $a\leq b$ and $a\in U$, then
$(\{b\},a)\in\mathcal D'_X$, because $a\in\da\{b\}$. Hence the
$\mathcal D'$-openness of $U$ implies $\{b\}\cap U\neq\emptyset$, so
$b\in U$.

Now fix $x\in U$. Put
\[
B_x=\{y\in X:x\in{\rm int}(\ua y)\}.
\]
By the property of continuous spaces, the set $B_x$ is directed, $B_x\leq x$, and
$B_x\to x$.  Hence $(B_x,x)\in\mathcal D'_X$. Since $U$ is
$\mathcal D'$-open and $x\in U$, there exists $y\in B_x\cap U$. As $U$ is
upper, $\ua y\subseteq U$; and as $y\in B_x$, we have
\[
x\in{\rm int}(\ua y)\subseteq \ua y\subseteq U.
\]
Thus every point of $U$ is an interior point, and $U$ is open. Hence
$X=\mathcal D'X=\mathbb D'X$ by Proposition \ref{RELA1}. Therefore every
continuous space belongs to $\bf Top_{\mathbb D'}$.

Since $\bf Top_{\mathbb D'}$ is coreflective, it is closed under colimits in
$\bf Top$. By Lemma \ref{COLI}, every continuous-space-generated space is a
colimit of continuous spaces, and therefore belongs to $\bf Top_{\mathbb D'}$.
Thus $\bf Top_{\mathsf{ConS}}\subseteq\bf Top_{\mathbb D'}$, and the reverse
inclusions above give the desired equality.
\end{proof}

\begin{theorem}\label{qcont-generated-dtop}

\( \bf Top_{\mathsf{QAlgS}} = \bf Top_{\mathsf{QConS}}=\bf Top_{\mathbb D}.
\)
\end{theorem}

\begin{proof}
By Proposition \ref{monotone core}, every one-point convergence space $(D\cup\{\infty\},\beta_D)$ is quasialgebraic. Thus $\mathbb D\subseteq\mathsf{QAlgS}\subseteq\mathsf{QConS}$, and therefore
\(
\bf Top_{\mathbb D}\subseteq \bf Top_{\mathsf{QAlgS}}\subseteq \bf Top_{\mathsf{QConS}}.
\)
Conversely, quasicontinuous spaces are monotone determined spaces. Hence $\bf Top_{\mathsf{QConS}} \subseteq \bf Top_{\mathbb D} $. 
The equality follows.
\end{proof}

\begin{remark}\label{continuous-qcontinuous-generated}
Combining Theorems \ref{cspaces-dprime-generated} and \ref{qcont-generated-dtop} with Proposition \ref{SUB}, one obtains the strict comparison
\[
\bf Top_{\mathsf{ConS}}
=\bf Top_{\mathbb D'}
\subset
\bf Top_{\mathsf{QConS}}
=\bf Top_{\mathbb D}.
\]
Thus continuous spaces and quasicontinuous spaces are both contained in monotone determined spaces, but they do not generate the same coreflective subcategory.
\end{remark}

For Scott spaces of dcpos, the situation is different: algebraic, continuous, quasicontinuous and arbitrary dcpos all generate the same coreflective subcategory of $\bf Top$.
\begin{theorem}\label{dcpo-scott-generated}
Let $\Sigma\mathsf{Dcpo}$ denote the class of all dcpos endowed with their Scott topologies. Let $\Sigma\mathsf{ADcpo}$, $\Sigma\mathsf{CDcpo}$ and $\Sigma\mathsf{QCDcpo}$ denote respectively the subclasses of algebraic dcpos, continuous dcpos and quasicontinuous dcpos, all with their Scott topologies. Then
\[
\bf Top_{\Sigma\mathsf{ADcpo}}
=\bf Top_{\Sigma\mathsf{CDcpo}}
=\bf Top_{\Sigma\mathsf{QCDcpo}}
=\bf Top_{\Sigma\mathsf{Dcpo}}.
\]
\end{theorem}

\begin{proof}
It suffices to show that every Scott space of a dcpo is a quotient of the
Scott space of an algebraic dcpo. Let $P$ be a dcpo, and let ${\rm Idl}(P)$ be
the dcpo of all non-empty directed lower subsets of $P$, ordered by inclusion.
This is an algebraic dcpo whose compact elements are the principal ideals
$\da x$.

Define
\[
q:{\rm Idl}(P)\longrightarrow P,\qquad q(I)=\bigvee I.
\]
The map $q$ is onto, since $q(\da x)=x$. We show that it is a quotient map for the Scott topologies. For $U\subseteq P$, if $U$ is Scott open, then $q^{-1}(U)$ is Scott open because $q$ preserves directed suprema.

Conversely, suppose that $q^{-1}(U)$ is Scott open in ${\rm Idl}(P)$. If $x\in U$ and $x\leq y$, then $\da x\subseteq\da y$; since $q^{-1}(U)$ is an upper set, $\da y\in q^{-1}(U)$, so $y\in U$. Thus $U$ is an upper set.

Let $D\subseteq P$ be directed and suppose that $\bigvee D\in U$. The family $(\da d)_{d\in D}$ is directed in ${\rm Idl}(P)$ and has supremum
\[
\bigcup_{d\in D}\da d.
\]
This ideal belongs to $q^{-1}(U)$, because its image under $q$ is $\bigvee D$. Since $q^{-1}(U)$ is Scott open, there exists $d\in D$ with $\da d\in q^{-1}(U)$, that is, $d\in U$. Therefore $U$ is inaccessible by directed suprema. Hence $U$ is Scott open.

We have proved
\[
U\in\sigma(P)\quad\Longleftrightarrow\quad q^{-1}(U)\in\sigma({\rm Idl}(P)).
\]
Thus $q:\Sigma{\rm Idl}(P)\to\Sigma P$ is a quotient map. Since generated categories are closed under quotient maps (see Lemma \ref{COLI}), every object of $\Sigma\mathsf{Dcpo}$ belongs to $\bf Top_{\Sigma\mathsf{ADcpo}}$, and the asserted equalities follow from the inclusions
\(
\Sigma\mathsf{ADcpo}\subseteq\Sigma\mathsf{CDcpo}
\subseteq\Sigma\mathsf{QCDcpo}\subseteq\Sigma\mathsf{Dcpo}.
\)
\end{proof}

\section{Relations among generated categories}

The inclusion relation of these categories of topological spaces is as follows.\ We use $\bf A \subseteq B$ to denote that $\bf A$ is a full subcategory of $\bf B$ and $\bf A \subset B$ to denote $\bf A \subseteq B$ and $\bf A \neq B$.\

\begin{remark} \label{TS}

Denote $\bf Seq$ the category of sequential spaces and $\bf TS$ the category of transfinite sequential spaces \cite{NRH}. For any directed poset $D$, let $\delta_D$ be the topology on $D\cup\{\infty\}$ generated by the subbase consisting of all $\ua x\cup\{\infty\}$ and $\{x\}$ for $x\in D$.\ Let 

$\mathbb{D}^{\prime \prime}  = \{(D \cup \{\infty\},\delta_{D}): D \in  \mathscr{D} \}$,\ $\mathcal{D}^{\prime \prime}_{X} = \mathcal{T}_{X}$ for any topological space $X$;

$\mathbb{I}^{\prime \prime}  = \{(I \cup \{\infty\},\delta_{I}): I \in  \mathscr{I} \}$,\ $\mathcal{I}^{\prime \prime}_{X} = \{(\{x_{i}\}_{i\in I},x):x\in X,\ \{x_i\} \text{ is a transfinite sequence},\ \{x_{i}\} \to x \}$;
 $\mathbb{N}^{\prime \prime}  = \{(N \cup \{\infty\},\delta_{N})\}$.\ 

Similarly,\ we have the conclusion that $\bf Top = \bf Top_{\mathbb{D}^{\prime \prime}} = \bf Top_{\mathcal{D}^{\prime \prime}}$,\ $\bf TS = \bf Top_{\mathbb{I}^{\prime \prime}} = \bf Top_{\mathcal{I}^{\prime \prime}}$,\ and $\bf Seq = \bf Top_{\mathbb{N}^{\prime \prime}} = \bf Top_{\mathcal{N}^{\prime \prime}}$.\
\end{remark}

The following proposition shows the inclusion relations among the categories of topological spaces generated by the six classes of one-point convergence spaces.

\begin{proposition} \label{SUB}
\begin{enumerate}
\item[(1)] $\bf Top_{\mathbb{N}} \subset \bf Top_{\mathbb{I} } \subset \bf Top_{\mathbb{D}}$,\ $\bf Top_{\mathbb{N}^{\prime}} \subset \bf Top_{\mathbb{I}^{\prime} } \subset \bf Top_{\mathbb{D}^{\prime}} $ 

\item[(2)] $\bf Top_{\mathbb{D}^{\prime}} \subset \bf Top_{\mathbb{D}} \subset \bf Top$,\ $\bf Top_{\mathbb{I}^{\prime}} \subset \bf Top_{\mathbb{I}} \subset \bf TS$,\ $\bf Top_{\mathbb{N}^{\prime}} \subset \bf Top_{\mathbb{N}} \subset \bf Seq$

\end{enumerate}

\end{proposition}

\begin{proof} 

(1) By Lemma \ref{COLI}, $\mathbb{N} \subset \mathbb{I} \subset \mathbb{D}$ implies $\bf Top_{\mathbb{N}} \subseteq \bf Top_{\mathbb{I}} \subseteq \bf Top_{\mathbb{D}}$.\ We give an example to show the strictness.\ Let $D=[\mathbb R]^{<\omega}$, ordered by inclusion.\ Considering $\mathbb{I}(D \cup \{\infty\},\beta_{D})$,\ we claim that $\{\infty\}$ is an open subset of it.\ Given any chain $I$ such that $I \subseteq D$,\ since any element $x \in I$ is a finite subset of $\mathbb R$,\ $I$ must be countable.\ Thus,\ there exists some $x \in D$ such that $x \notin \da I$.\ It follows that $I \cap (\ua x \cup \{\infty\}) = \emptyset$ and $I \not \ra \infty$.\ Therefore,\ $\{\infty\}$ is open in $\mathbb{I}(D \cup \{\infty\},\beta_{D}) $.\
Since $\{\infty\}$ is not open in $(D \cup \{\infty\},\beta_{D})$,\ then,\ $(D \cup \{\infty\},\beta_{D})$ is in $\bf Top_{\mathbb{D}}$ but not in $\bf Top_{\mathbb{I}}$.\ The proof for $\bf Top_{\mathbb{N}} \subset \bf Top_{\mathbb{I}}$ is similar,\ by replacing $D = \aleph_{1}$,\ where $\aleph_{1}$ is the least uncountable ordinal.\  

Similarly,\ $\mathbb{N}^{\prime} \subset \mathbb{I}^{\prime} \subset \mathbb{D}^{\prime} $ implies $\bf Top_{\mathbb{N}^{\prime}} \subseteq \bf Top_{\mathbb{I}^{\prime}} \subseteq \bf Top_{\mathbb{D}^{\prime}}$.\ $(D \cup \{\infty\},\gamma_{D})$, with $D=[\mathbb R]^{<\omega}$, is an example for $\bf Top_{\mathbb{I}^{\prime}} \neq \bf Top_{\mathbb{D}^{\prime}}$ and  $(\aleph_{1} \cup \{\infty\},\gamma_{\aleph_{1}})$ is an example for $\bf Top_{\mathbb{N}^{\prime}} \neq \bf Top_{\mathbb{I}^{\prime}}$.

(2) $\bf Top_{\mathbb{D}^{\prime}} \subseteq \bf Top_{\mathbb{D}} \subseteq \bf Top$ (resp.,\ $\bf Top_{\mathbb{I}^{\prime}} \subseteq \bf Top_{\mathbb{I}} \subseteq \bf TS$,\ $\bf Top_{\mathbb{N}^{\prime}} \subseteq \bf Top_{\mathbb{N}} \subseteq \bf Seq$) can be gained by combining Lemma \ref{CONV},\ Proposition \ref{RELA},\ Remark \ref{TS} and the fact that for any topological space $X$,\ $\mathcal{D}^{\prime}_{X} \subseteq  \mathcal{D}_{X} \subseteq \mathcal{D}^{\prime \prime}_{X}$ (resp.,
 $\mathcal{I}^{\prime}_{X} \subseteq  \mathcal{I}_{X} \subseteq \mathcal{I}^{\prime \prime}_{X}$,\ $\mathcal{N}^{\prime}_{X} \subseteq  \mathcal{N}_{X} \subseteq \mathcal{N}^{\prime \prime}_{X}$).\

Let $X = (N \cup \{\infty\},\beta_{N})$.\ Then $\{ \infty \}$ is open in $\mathbb{D}^{\prime} X$,\ but not open in $X$.\ Then $\bf Top_{\mathbb{D}^{\prime}} \subset \bf Top_{\mathbb{D}}$.\ It is also an example for $\bf Top_{\mathbb{I}^{\prime} } \subset \bf Top_{\mathbb{I} } $ and $\bf Top_{\mathbb{N}^{\prime} } \subset \bf Top_{\mathbb{N} }$.\ Similarly,\ let $Y = (N \cup \{\infty\},\delta_{N})$.\ It is also an example for $\bf Top_{\mathbb{D}} \subset \bf Top$,\ $\bf Top_{\mathbb{I} } \subset \bf TS  $ and $\bf Top_{\mathbb{N} } \subset \bf Seq $
\end{proof}

It is known that sequential spaces are compactly generated. Therefore, for countable cases, the above generated spaces are compactly generated. However, for uncountable cases, the generated spaces may not be compactly generated. The following theorem provides a counterexample already in the chain-generated category $\bf Top_{\mathbb{I}^{\prime}}$.

We first prove the following lemma, which is a closed-filtration version of the
ordinal Heine--Borel theorem for \(\omega_1\): compact subsets of
\(\omega_1\), with the order topology, are bounded.

\begin{lemma}\label{omega-one-closed-lemma}
Let $K$ be a compact Hausdorff space, and let $(C_\alpha)_{\alpha<\omega_1}$ be an increasing family of closed subsets of $K$ such that $C_\alpha \subseteq C_\beta$ for $\alpha < \beta$, and
\[
C_\lambda=\bigcup_{\alpha<\lambda}C_\alpha
\]
for every countable limit ordinal $\lambda<\omega_1$. Then $\bigcup_{\alpha<\omega_1}C_\alpha$ is closed in $K$.
\end{lemma}

\begin{proof}
Let $U=\bigcup_{\alpha<\omega_1}C_\alpha$ and suppose, for a contradiction, that $x\in\overline U\setminus U$. 
We recursively construct a strictly increasing sequence
\(
\alpha_0<\alpha_1<\alpha_2<\cdots
\)
and a decreasing sequence of open neighbourhoods of \(x\),
\(
V_0\supseteq V_1\supseteq V_2\supseteq\cdots,
\)
such that
\[
\overline{V_{n+1}}\subseteq V_n,
\overline{V_n}\cap C_{\alpha_n}=\emptyset,
V_n\cap C_{\alpha_{n+1}}\neq\emptyset.
\]

Choose any \(\alpha_0<\omega_1\). Since \(x\notin C_{\alpha_0}\) and
\(C_{\alpha_0}\) is closed, regularity of $K$ gives an open neighbourhood \(V_0\) of
\(x\) such that
\(
\overline{V_0}\cap C_{\alpha_0}=\emptyset.
\)
Since \(x\in\overline{U}\), we have
\(
V_0\cap U\neq\emptyset.
\)
Since \(V_0\cap C_{\alpha_0}=\emptyset\) and the family ${C_\alpha}_{\alpha<\omega_1}$ is increasing, there
exists
\(
\alpha_1>\alpha_0
\)
such that
\(
V_0\cap C_{\alpha_1}\neq\emptyset.
\)

Now, by induction, suppose that \(\alpha_n\) and \(V_n\) have been chosen. Since
\(x\notin C_{\alpha_{n+1}}\) and \(C_{\alpha_{n+1}}\) is closed, we can choose an open neighbourhood \(V_{n+1}\) of \(x\) such that
\(
\overline{V_{n+1}}\subseteq V_n
\)
and
\(
\overline{V_{n+1}}\cap C_{\alpha_{n+1}}=\emptyset.
\)
Again, since \(x\in\overline U\), we have
\(
V_{n+1}\cap U\neq\emptyset.
\)
Since \(V_{n+1}\cap C_{\alpha_{n+1}}=\emptyset\) and the family is increasing, there exists
\(
\alpha_{n+2}>\alpha_{n+1}
\)
such that
\(
V_{n+1}\cap C_{\alpha_{n+2}}\neq\emptyset.
\)
This completes the construction.

Let
\(
\lambda=\sup_{n<\omega}\alpha_n.
\)
Then \(\lambda<\omega_1\), and \(\lambda\) is a countable limit ordinal. By the assumption,
since \((\alpha_n)\) is cofinal in \(\lambda\), we have
\[
C_\lambda=\bigcup_{n<\omega}C_{\alpha_n}.
\]

For each \(n\), choose
\(
y_n\in V_n\cap C_{\alpha_{n+1}}.
\)
If \(m\geq n\), then \(V_m\subseteq V_n\), and therefore
\(
y_m\in V_n.
\)
Also \(y_m\in C_\lambda\). Hence for each \(n\),
\(
C_\lambda\cap\overline{V_n}\neq\emptyset.
\)
The sets
\(
C_\lambda\cap\overline{V_n}
\)
are nonempty closed subsets of the compact space \(C_\lambda\), and they form a
decreasing sequence. Thus
\[
\bigcap_{n<\omega}\bigl(C_\lambda\cap\overline{V_n}\bigr)\neq\emptyset.
\]
Choose
\[
y\in \bigcap_{n<\omega}\bigl(C_\lambda\cap\overline{V_n}\bigr).
\]
Since
\[
C_\lambda=\bigcup_{n<\omega}C_{\alpha_n},
\]
there exists \(N<\omega\) such that
\(
y\in C_{\alpha_N}.
\)
But also
\(
y\in\overline{V_N},
\)
while
\(
\overline{V_N}\cap C_{\alpha_N}=\emptyset,
\)
a contradiction. Therefore \(U\) is closed.
\end{proof}

\begin{theorem}\label{omega-one-common-counterexample}
Let
\[
Y_{\omega_1}=\omega_1\cup\{\infty\}
\]
carry the topology generated by the sets
\[
\ua\alpha\cup\{\infty\}
=\{\beta<\omega_1:\alpha\leq \beta\}\cup\{\infty\},
\qquad \alpha<\omega_1.
\]
Then
\[
Y_{\omega_1}\in \bf Top_{\mathbb I'}
\quad\text{but}\quad
Y_{\omega_1}\notin \bf Top_{\mathbb K}.
\]
\end{theorem}

\begin{proof}
The space $Y_{\omega_1}$ is exactly the strong one-point convergence space
$(\omega_1\cup\{\infty\},\gamma_{\omega_1})$ associated with the chain
$\omega_1$. Hence it is a generator in $\mathbb I'$ and therefore belongs to
$\bf Top_{\mathbb I'}$.

The set $\{\infty\}$ is not open in $Y_{\omega_1}$, because every non-empty basic open set (and hence every neighbourhood) containing $\infty$ must also contain some ordinal $\alpha<\omega_1$.

To show $Y_{\omega_1}\notin \bf Top_{\mathbb K}$, we demonstrate that $\{\infty\}$ is $K$-open but not open. Thus it is enough to show that, for every continuous map $p:K\to Y_{\omega_1}$ from a compact Hausdorff space $K$, the preimage $p^{-1}(\{\infty\})$ is open in $K$. 

Let $p:K\to Y_{\omega_1}$ be any continuous map from a compact Hausdorff space, and put $B=p^{-1}(\{\infty\})$. For each $\alpha<\omega_1$, define the open set
\[
W_\alpha=p^{-1}(\ua\alpha\cup\{\infty\}),
\]
and its closed complement $C_\alpha=K\setminus W_\alpha$. 
Clearly, $(C_\alpha)_{\alpha<\omega_1}$ is an increasing family of closed sets in $K$. If $\lambda<\omega_1$ is a countable limit ordinal, then for every $k\in K$ one has
\[
k\in C_\lambda \Longleftrightarrow p(k)\notin \ua\lambda\cup\{\infty\} \Longleftrightarrow p(k)=\gamma<\lambda.
\]
Since $\lambda$ is a limit ordinal, $p(k)=\gamma<\lambda$ implies $p(k) < \gamma+1 \leq \lambda$, meaning $k \in C_{\gamma+1}$. Thus,
\[
p(k)<\lambda \Longleftrightarrow k\in \bigcup_{\alpha<\lambda}C_\alpha.
\]
This verifies the hypothesis of Lemma \ref{omega-one-closed-lemma}, ensuring that $U=\bigcup_{\alpha<\omega_1}C_\alpha$ is closed in $K$. 

For any $k \in K$, if $p(k)=\infty$, then $p(k) \in \ua\alpha\cup\{\infty\}$ for all $\alpha$, so $k\notin C_\alpha$ for all $\alpha$. Conversely, if $p(k)=\gamma<\omega_1$, then $k\in C_{\gamma+1} \subseteq U$. Therefore, $U=K\setminus B$. Since $U$ is closed, its complement $B = p^{-1}(\{\infty\})$ is open in $K$.

Because $p^{-1}(\{\infty\})$ is open in $K$ for every continuous map $p$ from a compact Hausdorff space, but $\{\infty\}$ is not open in $Y_{\omega_1}$, it follows by definition that $Y_{\omega_1}\notin \bf Top_{\mathbb K}$.
\end{proof}

Theorem \ref{omega-one-common-counterexample} gives a counterexample already generated by chains. Hence
\[
\bf Top_{\mathbb I'}\not\subseteq \bf Top_{\mathbb K}.
\]
Using the inclusions in Proposition \ref{SUB}, the same example also yields
\[ 
\bf Top_{\mathbb D'}\not\subseteq \bf Top_{\mathbb K}
\qquad\text{and}\qquad
\bf Top_{\mathbb I}\not\subseteq \bf Top_{\mathbb K}.
\]
Thus the compact Hausdorff generated category does not contain the chain-determined or monotone determined categories considered here. On the other hand, the present example is a $\gamma$-type one-point convergence space, hence a continuous space by Proposition \ref{gamma-cspaces}; therefore the class of continuous spaces, and the category $\bf Top_{\mathsf{ConS}}=\bf Top_{\mathbb D'}$, are not contained in $\bf Top_{\mathbb K}$.

The second consequence clarifies why the chain-generated categories are genuinely different from the directed ones outside the dcpo setting. It is well known that, for posets, directed completeness and chain completeness coincide \cite{GM1976}. Given a poset $P$, a subset $U$ of $P$ is called chain open if for any chain $I$ such that $\sup I = x$ and $x \in U$, one has $I \cap U \neq \emptyset$.  
For any dcpo $L$, the chain-open subsets of $L$ agree with the Scott-open subsets. Equivalently,
\[
\mathbb{I}(L,\upsilon(L))=\mathbb{D}(L,\upsilon(L)).
\]
For general posets, however, the chain-open and Scott-open subsets need not coincide. The following example isolates this separation.

\begin{example}

Let $D=[\mathbb R]^{<\omega}$ be the set of all finite subsets of $\mathbb R$, ordered by inclusion, and let $\mathcal{M}$ be the set of all maximal chains of $D$.\ For each $I \in \mathcal{M}$,\ we define an element $\alpha_{I}$ that is different from $\mathbb R$ and all elements in $D$.\ Let $S = \{\mathbb R\} \cup D \cup \{\alpha_{I}\}_{I \in \mathcal{M}}$.\
We define the order $\leq$ on $S$ as follows:
\begin{enumerate}
\item For $x,y \in D \cup\{ \mathbb R \}$:\ the order is the inclusion order of sets,\ i.e.,\ $x \leq y$ iff $x \subseteq y$;

\item For $x \in D$ and $y = \alpha_{I}$:
$x \leq y$ iff $\exists z \in I$, $x \leq z $;

\item For $x,y \in \{\mathbb R\} \cup \{\alpha_{I}\}_{I \in \mathcal{M}}$: $x \leq y$ iff $x = y$.
\end{enumerate}
It is easy to check that $(S,\leq)$ is a well defined poset.\ We show that $D$ is a chain closed subset,\ but not a Scott closed subset.\ Finite chains in $D$ have their suprema in $D$. If $C$ is an infinite chain in $D$,\ extend $C$ to a maximal chain $I\in\mathcal M$. Then both $\mathbb R$ and $\alpha_I$ are upper bounds of $C$, and they are incomparable. Hence $C$ has no supremum in $S\setminus D$, so $D$ is chain closed.\ For any $\alpha_{I}$,\ since $\da \alpha_{I}$ is countable,\ there exists an element $y$ in $D$ such that $y \nleq \alpha_{I}$,\ i.e.,\ $\alpha_{I}$ is not an upper bound of $D$.\ Hence,\ $\mathbb R$ is the supremum of $D$ and then $D$ is not Scott closed.\
In fact,\ every $D \cup A$,\ where $A$ is a subset of $\{\alpha_{I}\}_{I \in \mathcal{M}}$,\ is chain closed but not Scott closed. $\qedsymbol$
\end{example}

The preceding example shows that the agreement over dcpos does not extend to arbitrary posets. More generally, for every dcpo $L$ one has
\[
\mathbb{I}^{\prime}(L,\upsilon(L))
=\mathbb{I}(L,\upsilon(L))
=\mathbb{D}(L,\upsilon(L))
=\mathbb{D}^{\prime}(L,\upsilon(L))
=(L,\sigma(L)).
\]

We finally record that the compact-generation theorem for Scott spaces of
dcpos cannot be extended from dcpos to arbitrary posets.

\begin{theorem}\label{scott-poset-not-kgenerated}
There exists a poset $P$ such that $(P,\sigma(P))$ is not
compact Hausdorff generated.
\end{theorem}

\begin{proof}
Let
\[
L=\omega_1\times\mathbb Q
\]
be ordered lexicographically:
\[
(\alpha,q)\leq(\beta,r)
\quad\Longleftrightarrow\quad
\alpha<\beta\quad\text{or}\quad(\alpha=\beta\text{ and }q\leq r).
\]
Adjoin a largest element $\infty$ and put $P=L\cup\{\infty\}$. This poset is
not a dcpo; for instance,
\[
\{(0,q):q\in\mathbb Q,\ q^2<2\}
\]
is directed but has no supremum in $P$.

For each $\alpha<\omega_1$, set
\[
U_\alpha=\{(\beta,q)\in\omega_1\times\mathbb Q:\alpha\leq\beta\}\cup\{\infty\}.
\]
We first show that $U_\alpha$ is Scott open. It is clearly an upper set. Let
$A\subseteq P$ be directed, suppose that $\sup A=s\in U_\alpha$, and assume
that $A\cap U_\alpha=\emptyset$. Then every element of $A$ lies in a block
strictly below $\alpha$. If $s=(\alpha,q)$, choose $q'<q$; then
$(\alpha,q')$ is still an upper bound of $A$, contradicting $s=\sup A$. If
$s=(\beta,q)$ with $\beta>\alpha$, then $(\alpha,0)$ is an upper bound of $A$
strictly below $s$. If $s=\infty$, the same element $(\alpha,0)$ is an upper
bound of $A$ strictly below $\infty$. These contradictions show that
$A\cap U_\alpha\neq\emptyset$. Hence $U_\alpha$ is Scott open.

The singleton $\{\infty\}$ is not Scott open: the set $L$ is directed,
$\sup L=\infty$, and $L\cap\{\infty\}=\emptyset$.

It remains to show that $\{\infty\}$ is open with respect to compact
Hausdorff probes. Let $p:K\to(P,\sigma(P))$ be continuous, where $K$ is
compact Hausdorff, and put $B=p^{-1}(\{\infty\})$. For each
$\alpha<\omega_1$, define
\[
C_\alpha=K\setminus p^{-1}(U_\alpha)
=p^{-1}\bigl(\{(\beta,q):\beta<\alpha\}\bigr).
\]
Each $C_\alpha$ is closed in $K$, and $(C_\alpha)_{\alpha<\omega_1}$ is
increasing. If $\lambda<\omega_1$ is a countable limit ordinal, then
\[
C_\lambda=\bigcup_{\alpha<\lambda}C_\alpha.
\]
Indeed, if $p(k)=\infty$, then $k$ belongs to no $C_\alpha$; if
$p(k)=(\gamma,q)$, then
\[
k\in C_\lambda
\Longleftrightarrow
\gamma<\lambda
\Longleftrightarrow
\exists\alpha<\lambda\;(\gamma<\alpha)
\Longleftrightarrow
k\in\bigcup_{\alpha<\lambda}C_\alpha.
\]
By Lemma \ref{omega-one-closed-lemma},
\[
C=\bigcup_{\alpha<\omega_1}C_\alpha
\]
is closed in $K$. Moreover, $C=K\setminus B$: if $p(k)=\infty$, then
$k\notin C$, while if $p(k)=(\gamma,q)\in L$, then $k\in C_{\gamma+1}$.
Thus $B$ is open in $K$.

Therefore $p^{-1}(\{\infty\})$ is open for every compact Hausdorff probe
$p:K\to(P,\sigma(P))$, but $\{\infty\}$ is not Scott open. Hence
$(P,\sigma(P))$ is not compact Hausdorff generated.
\end{proof}

We now summarize where posets endowed with the Scott topology and Scott spaces
of dcpos sit among the generated categories. Let $\Sigma\mathsf{Poset}$ denote
the class of all posets endowed with their Scott topologies, and let
$\Sigma\mathsf{Dcpo}$ denote the subclass of Scott spaces of dcpos.

First, every poset endowed with the Scott topology is strongly directed
determined. Indeed, let $P$ be a poset and consider $\Sigma P=(P,\sigma(P))$. If
$U\subseteq P$ is $\mathcal D'$-open, then $U$ is an upper set: whenever
$x\in U$ and $x\leq y$, the singleton pair $(\{y\},x)$ belongs to
$\mathcal D'_{\Sigma P}$, and hence $y\in U$. If $D\subseteq P$ is directed,
$\bigvee D=x$ exists and $x\in U$, then $D\to x$ in $\Sigma P$ and
$D\leq x$, so $(D,x)\in\mathcal D'_{\Sigma P}$; hence
$D\cap U\neq\emptyset$. Thus $U$ is Scott open. Conversely, every Scott open
set is clearly $\mathcal D'$-open, by upperness and Scott convergence of
directed suprema. Therefore
\[
\Sigma\mathsf{Poset}\subseteq \bf Top_{\mathbb D'}
\subseteq \bf Top_{\mathbb D}.
\]

\begin{remark}
This inclusion cannot in general be improved to the chain-generated category
$\bf Top_{\mathbb I}$. For the poset $S$ constructed above from
$D=[\mathbb R]^{<\omega}$ and its maximal chains, put
\[
U=S\setminus D=\{\mathbb R\}\cup\{\alpha_I:I\in\mathcal M\}.
\]
Then $U$ is not Scott open, since $D$ is directed, $\sup D=\mathbb R\in U$,
and $D\cap U=\emptyset$. We claim that $U$ is nevertheless
$\mathcal I$-open in $\Sigma S$. Let $C$ be a chain in $S$ such that
$C\to x\in U$. If $C$ is not eventually in $U$, then, since $U$ is an upper
set, $C\subseteq D$. Extend $C$ to a maximal chain $I\in\mathcal M$. Then
$C\subseteq\da\alpha_I$. If $x=\mathbb R$ or $x=\alpha_J$ with $J\neq I$, the
Scott open set $S\setminus\da\alpha_I$ contains $x$ and is disjoint from $C$.
If $x=\alpha_I$, the Scott open set $S\setminus\da\mathbb R$ contains $x$ and
is disjoint from $C$. In each case we contradict $C\to x$. Hence $C$ is
eventually in $U$, and $U$ is $\mathcal I$-open. Consequently
\[
\Sigma\mathsf{Poset}\not\subseteq\bf Top_{\mathbb I},
\]
and since $\bf Top_{\mathbb I'}\subseteq\bf Top_{\mathbb I}$, it follows also
that $\Sigma\mathsf{Poset}\not\subseteq\bf Top_{\mathbb I'}$.

For dcpos the situation is stronger. If $L$ is a dcpo, then the
chain-open subsets of $L$ coincide with the Scott-open subsets
\cite{GM1976}. The same argument as above, with chains in place of arbitrary
directed subsets, gives
\[
\mathbb I'(\Sigma L)=\mathbb I(\Sigma L)
=\mathbb D'(\Sigma L)=\mathbb D(\Sigma L)=\Sigma L.
\]
Therefore
\[
\Sigma\mathsf{Dcpo}\subseteq \bf Top_{\mathbb I'}.
\]
\end{remark}

Together with the compact-generation theorem for Scott spaces of dcpos
\cite[Theorem 4.7]{CG2004}, this gives
\[
\Sigma\mathsf{Dcpo}\subseteq
\bf Top_{\mathbb K}\cap \bf Top_{\mathbb I'}.
\]
By contrast, 
\(
\Sigma\mathsf{Poset}\not\subseteq \bf Top_{\mathbb K} \ \text{and} \  \Sigma\mathsf{Poset}\not\subseteq \bf Top_{\mathbb I}.
\)

We end by summarizing the inclusion relations among the categories discussed above. An arrow $\bf A\to\bf B$ means that $\bf B$ is a full subcategory of $\bf A$.

\begin{center}
\begin{tikzpicture}[
  node distance=1.35cm and 1.9cm,
  cat/.style={draw, rounded corners=2pt, align=center, inner xsep=7pt, inner ysep=4pt, font=\small, fill=gray!6},
  main/.style={cat, fill=blue!6},
  sub/.style={cat, fill=green!6},
  poscat/.style={cat, fill=orange!8},
  arr/.style={->, >=stealth, line width=.45pt}
]

\node[main] (top) at (4,4) {$\bf Top$};
\node[main] (topk) at (0,3) {$\bf Top_{\mathbb K}$};
\node[main] (ts) at (2,3) {$\bf TS$};
\node[main] (seq) at (0,2) {$\bf Seq$};
\node[sub] (topd) at (4,3) {$\bf Top_{\mathbb D}$};
\node[sub] (topdp) at (4,2) {$\bf Top_{\mathbb D'}$};
\node[sub] (topi) at (2,2) {$\bf Top_{\mathbb I}$};
\node[sub] (topip) at (2,1) {$\bf Top_{\mathbb I'}$};
\node[sub] (topn) at (0,1) {$\bf Top_{\mathbb N}$};
\node[sub] (topnp) at (0,0) {$\bf Top_{\mathbb N'}$};

\node[poscat] (poss) at (6,2) {$\Sigma\mathsf{Poset}$};
\node[poscat] (dcpo) at (6,1) {$\Sigma\mathsf{Dcpo}$};

\draw[arr] (top) -- (ts);
\draw[arr] (top) -- (topk);
\draw[arr] (topk) -- (seq);
\draw[arr] (ts) -- (seq);
\draw[arr] (seq) -- (topn);
\draw[arr] (topn) -- (topnp);
\draw[arr] (ts) -- (topi);
\draw[arr] (topi) -- (topip);
\draw[arr] (topi) -- (topn);
\draw[arr] (topip) -- (topnp);
\draw[arr] (top) -- (topd);
\draw[arr] (topd) -- (topdp);
\draw[arr] (topd) -- (topi);
\draw[arr] (topdp) -- (topip);
\draw[arr] (topdp) -- (poss);
\draw[arr] (topip) to[bend right=8] (dcpo);
\draw[arr] (topk.south west) .. controls (-1.0,-0.8) and (5.5,-0.8) .. (dcpo.south west);
\draw[arr] (poss) -- (dcpo);

\end{tikzpicture}
\end{center}

\end{document}